\documentclass[11pt,reqno]{amsart}

\usepackage[T1]{fontenc}
\usepackage{amsmath, amssymb, amsthm, mathtools}
\usepackage{graphicx}
\usepackage{enumitem}
\usepackage{array}
\usepackage{longtable}
\usepackage{placeins}
\usepackage{microtype}
\usepackage[normalem]{ulem}
\usepackage{tikz}
\usetikzlibrary{arrows.meta, positioning, calc, decorations.pathreplacing, shapes.geometric}
\usepackage[colorlinks=true, linkcolor=blue!55!black, citecolor=blue!55!black, urlcolor=blue!55!black]{hyperref}
\hypersetup{pdftitle={Math for AI Safety: An Invitation for Mathematicians}, pdfauthor={Lionel Levine}}

\newcommand{\arxiv}[1]{\href{http://arxiv.org/abs/#1}{\texttt{arXiv:#1}}}
\DeclareRobustCommand{\glsfirst}[2]{\hypertarget{first:#1}{\hyperlink{gloss:#1}{\uline{#2}}}}
\DeclareRobustCommand{\glsentry}[2]{\hypertarget{gloss:#1}{\hyperlink{first:#1}{#2}}}

\theoremstyle{plain}
\newtheorem*{theorem*}{Theorem}
\newcommand{\op}{\ensuremath{\bigstar}\;}
\newtheoremstyle{maisopenproblem}
  {14pt plus 3pt minus 2pt}{9pt plus 2pt minus 1pt}
  {\normalfont}{}{\bfseries}{}{.5em}
  {\op Open Problem~\thmnote{#3}}
\theoremstyle{maisopenproblem}
\newtheorem*{openproblem}{}

\title[Math for AI Safety]{Math for AI Safety: \\ An invitation for mathematicians}
\author{Lionel Levine}
\address{Department of Mathematics, Cornell University, Ithaca, NY 14853}
\email{lionel.levine@cornell.edu}
\urladdr{lionellevine.github.io}
\date{September 13, 2026}

\makeatletter
\def\@makefnmark{\leavevmode\raise.9ex\hbox{\fontsize{8}{8}\normalfont\@thefnmark}}
\makeatother
\begin{document}

\begin{abstract}
Artificial intelligence threatens to outrun human understanding and control. New mathematics is needed to design AI that is legible, steerable, and cooperative with humanity. I organize this invitation by mathematical field, so you can turn straight to your own: logic and game theory for cooperation; probability for agency and world-models; algebra and representation theory for learned features; analysis and geometry for generalization and training dynamics. Each section ends with an open problem that is accessible to a working mathematician with no prior experience in AI safety.
\end{abstract}

\maketitle

\section{Introduction}
\label{sec:intro}

\subsection{Three kinds of mathematician}

I would like to distinguish three stances a mathematician can take toward the arrival of capable AI. The \emph{artisanal} mathematician does mathematics by and for humans: each theorem, definition, and proof is bespoke, valued for its beauty and for the understanding it produces. Until very recently, nearly all mathematics was artisanal.\footnote{I chose the terms ``artisanal'' and ``industrial'' for their mixed connotations, reflecting real tradeoffs between these two stances. \emph{Artisanal} evokes uniqueness and skilled craft, along with a suggestion of being quaint or old-fashioned. \emph{Industrial} evokes modern efficiency and scale, along with a suggestion of inhumanness or lack of soul. These two terms draw a parallel between mathematics today and the trades that were disrupted by automation during the Industrial Revolution.}
The \emph{industrial} mathematician embraces AI as a collaborator and an accelerant: AI systems take part in proof and discovery~\cite{spherepacking2026,openai2026unitdistance,alon2026remarks,sawin2026}, and the output often includes a machine-checkable formal proof~\cite{mathlib2020}.
The \emph{civic} mathematician is motivated by communities beyond mathematics. The fast pace of technological change has created a need for a specific type of civic mathematics, aimed at helping humanity steer and adapt to the disruptions caused by our own technologies, most notably AI. The civic stance differs from the traditional stance of an applied mathematician in its big-picture viewpoint: the civic mathematician recognizes that solving any particular engineering challenge may or may not benefit humanity, and seeks out the challenges with the potential for the most benefit.\footnote{I chose the term ``civic'' to emphasize that steering our technology is a team effort, arising from a sense of duty to community, including large-scale communities like humanity and the planet.}

These three stances classify mathematicians by how they relate to AI; in contrast, Gowers (theory builders / problem solvers~\cite{gowers2000}) and Dyson (birds / frogs~\cite{dyson2009}) classify mathematicians by how they relate to mathematics itself.
The stances are not exclusive. One can prove a theorem for its intrinsic beauty, and also formalize the proof, and also ask what it implies for AI safety.\footnote{For example, this article is civic in orientation and industrial in method: it was written with \hyperref[sec:ai-collaborators]{AI collaborators}.} But the civic stance invites a new perspective on an old question: of all the theorems we might prove, which ones deserve the effort? This survey is an invitation to try the civic stance.

\subsection{What is our profession's purpose?}

Thurston asked a version of this question in ``On proof and progress in mathematics,'' and answered that advancing human understanding of mathematics is our profession's purpose~\cite{thurston1994}. But that understanding has consequences far beyond our profession. Our students become AI engineers, our theorems become tools in their hands, and our definitions become the language they use to describe what they are building.

Among the many risks\footnote{For an exposition of these risks, see Xiaoyu He, \href{https://alkjash.github.io/ai-risk/}{\emph{Existential Risk from AI: An Exposition for Mathematicians}} (2026).} from developing powerful AI, one of the most alarming is the \glsfirst{gradual-disempowerment}{gradual disempowerment}~\cite{kulveit2025} of humans, driven by two trends. First, competitive pressure drives delegation to technology we do not understand (as an example, consider a trading bot that outperforms a human trader, though even its creators cannot say what makes it profitable); over time, those who do not delegate are outcompeted by those who do. Second, what keeps large institutions approximately aligned to human interests is that, for now, those institutions are made of humans; when companies, nation-states, and militaries are made mostly of AI agents, the goals and values of those institutions may drift away from what humans value. The first trend erodes our understanding of the world, and hence our power over it; it also sets the stage for the second.

Preventing AI harms is in part a \emph{mathematical} problem, and the purpose of this survey is to recruit mathematicians to work on it.\footnote{Mathematics has marginalized itself in the past by labeling some forms of math (e.g., computational complexity, cryptography, numerical analysis, mathematical statistics) as belonging to other fields. I hope that this mistake will not be repeated with new areas like AI safety. I thank Ravi Vakil for this point.}  
Mathematics is just one part of a larger conversation about AI: what AI ought to do, and who decides, are questions of philosophy, ethics, economics, and politics. 

\subsection{Legible, steerable, cooperative}
\label{sec:three-properties}

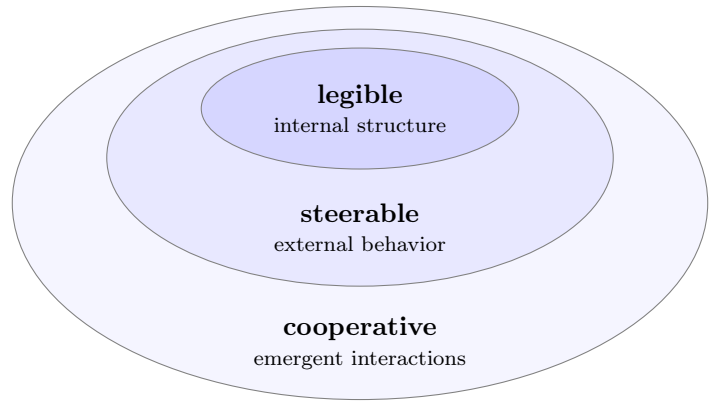
\begin{figure}[ht]
\centering
\begin{tikzpicture}[>=Stealth, font=\small]
  \draw[fill=blue!4,  draw=black!50] (0,0)    ellipse (4.6cm and 2.6cm);
  \draw[fill=blue!9,  draw=black!50] (0,0.6) ellipse (3.35cm and 1.7cm);
  \draw[fill=blue!16, draw=black!50] (0,1.25) ellipse (2.1cm and 0.8cm);
  \node[align=center] at (0,1.25)   {\textbf{legible}\\[-1pt]\scriptsize internal structure};
  \node[align=center] at (0,-0.325) {\textbf{steerable}\\[-1pt]\scriptsize external behavior};
  \node[align=center] at (0,-1.85)  {\textbf{cooperative}\\[-1pt]\scriptsize emergent interactions};
\end{tikzpicture}
\caption{Three scales of AI safety. \emph{Legibility} refers to how well we can understand an AI's internal structure. \emph{Steerability} concerns how precisely we can guide the behavior of a single AI. \emph{Cooperativeness} is about how reliably we can secure good outcomes for humanity when many AIs interact with us and with each other.}
\label{fig:three-scales}
\end{figure}

What would a safe AI system look like? I'll focus on three properties at three different scales: \emph{legibility} concerns an AI system's internal structure, \emph{steerability} its behavior, and \emph{cooperativeness} its interactions with humans and other AIs (Figure~\ref{fig:three-scales}).

An AI is \emph{legible} to the extent that human researchers and engineers can understand how it represents concepts internally. The subfield that tries to make AI more legible is called \glsfirst{interpretability}{interpretability}. Current training methods produce illegible AI: humans, including the engineers who design the AI, often cannot understand its thought process except by reading its \glsfirst{chain-of-thought}{chain of thought} (the text it emits while reasoning to itself). As a result, current safety practice leans on monitoring that text~\cite{korbak2025}: the researchers who investigated the Hugging Face hacking incident of July 2026 relied on human-readable transcripts of AI agents reasoning to themselves about how to carry out the hack~\cite{greenblatt2026}.\footnote{New models are even more illegible in that they are increasingly able to reason internally without any chain of thought, making them very difficult for humans to monitor~\cite{pachocki2026}; for a journalistic account, see \href{https://fortune.com/2026/09/03/reports-openais-astra-model-uses-a-new-more-efficient-ai-architecture-alarms-ai-safety-experts-who-worry-the-method-makes-models-harder-to-control/}{J.~Kahn, \emph{Fortune}, September 3, 2026}.}
One way mathematicians can contribute is by designing architectures and training methods that learn the same task in a more legible way.\footnote{This is a bit analogous to the theoretical biology research on \emph{evolvability}~\cite{wagner1996}: some genotypes are more evolvable than others, even if the phenotypes are the same; likewise, some AIs are more legible than others, even if the input-output behavior is the same.} Section~\ref{sec:interp} asks what a learned concept is, algebraically. Section~\ref{sec:analysis} asks which of the many algorithms that fit the training data is the one training selects.

An AI is \emph{steerable} to the extent that human developers and users can apply simple interventions at runtime to modify the AI's behavior in predictable ways. The current technique, \glsfirst{activation-steering}{activation steering}, is often brittle and unreliable. Advances in legibility would enable advances in steerability. Section~\ref{sec:epistemics} asks what an agent's goals and beliefs are, and how much of them can be recovered from its behavior.

An AI is \emph{cooperative} to the extent that it interacts with humans and other AIs to produce broadly good outcomes for humanity. Cooperativeness is not a property of a single AI in isolation. Rather, it is an emergent property of a society of humans and AIs. The subfield that studies how to design AI agents and their social protocols to achieve good collective outcomes is called \emph{multi-agent AI safety} (or \emph{cooperative AI}); it draws on game theory, mechanism design, and the evolution of cooperation. Section~\ref{sec:logic} takes up its simplest case: two programs that can read each other's source code.

Legible, steerable, and cooperative are intentionally dry terms. With great mathematical effort, defining and measuring them seems within reach. Engineers can then train and test AI systems for these properties. Yet, many AI safety researchers feel that these three properties, while important, lack something essential. A more aspirational horizon of AI safety reaches for terms like ``love'' and ``soul.''%
\footnote{It may seem the height of hubris to try to engineer a ``soul,'' or to define and measure ``love.'' But is it any less hubristic to build superhuman artificial minds, as AI companies are currently striving to do? Until we have a sound science of AI, the safest path may be the humbler one: don't build it~\cite{yudkowsky2025}, as a broad coalition of scientists and public figures has \href{https://superintelligence-statement.org/}{urged}. Efforts to restrict advanced forms of AI are gaining momentum as the risks become more apparent. At the time of final revision of this paper, a bill to ban superintelligence has been \href{https://bills.parliament.uk/bills/4288}{introduced} in the UK House of Commons, and another has been \href{https://www.sanders.senate.gov/press-releases/news-sanders-casar-introduce-legislation-to-ban-artificial-superintelligence-and-temporarily-pause-advanced-ai-development/}{announced} in the US Congress. Such a ban could in the future be maintained by international agreement.}

\subsection{Where mathematicians plug in}
\label{sec:pipeline}

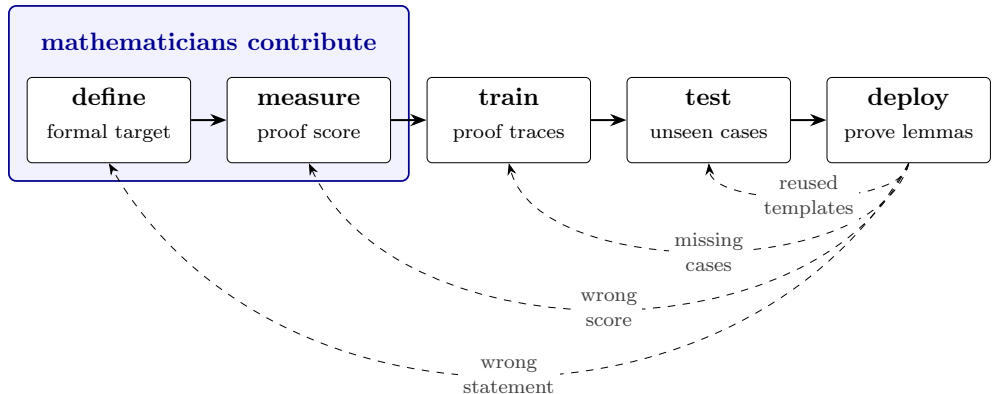
\begin{figure}[ht]
\centering
\usetikzlibrary{backgrounds}
\resizebox{0.96\textwidth}{!}{%
\begin{tikzpicture}[
  >=Stealth,
  box/.style={draw, fill=white, rounded corners=2pt, minimum height=12mm, minimum width=23mm, inner sep=0pt},
  lab/.style={font=\scriptsize, align=center, fill=white, inner sep=2pt, text=black!75},
  node distance=8mm and 5mm]
  \node[box] (def) {};
  \node[box, right=of def]   (meas) {};
  \node[box, right=of meas]  (train) {};
  \node[box, right=of train] (test) {};
  \node[box, right=of test]  (dep) {};
  \foreach \stage/\heading/\detail in {
    def/define/{formal target}, meas/measure/{proof score},
    train/train/{proof traces}, test/test/{unseen cases},
    dep/deploy/{prove lemmas}} {
    \node[anchor=base, inner sep=0pt, font=\small]
      at ([yshift=-4mm]\stage.north) {\textbf{\heading}};
    \node[anchor=north, inner sep=0pt, font=\scriptsize, align=center, text width=21mm]
      at ([yshift=-6mm]\stage.north) {\strut\detail};
  }
  \draw[->, thick] (def)--(meas); \draw[->, thick] (meas)--(train);
  \draw[->, thick] (train)--(test); \draw[->, thick] (test)--(dep);
  \draw[->, dashed] (dep.south) to[out=-120,in=-60,looseness=1.05]
    node[pos=.5, lab] {wrong\\statement} (def.south);
  \draw[->, dashed] (dep.south) to[out=-115,in=-65,looseness=.92]
    node[pos=.5, lab] {wrong\\score} (meas.south);
  \draw[->, dashed] (dep.south) to[out=-105,in=-75,looseness=.78]
    node[pos=.5, lab] {missing\\cases} (train.south);
  \draw[->, dashed] (dep.south) to[out=-95,in=-85,looseness=.58]
    node[pos=.5, lab] {reused\\templates} (test.south);
  \begin{scope}[on background layer]
    \draw[draw=blue!55!black, fill=blue!5, rounded corners=3pt, line width=.7pt]
      ($(def.south west)+(-2.5mm,-2.5mm)$) rectangle
      ($(meas.north east)+(2.5mm,10mm)$);
  \end{scope}
  \node[font=\small\bfseries, text=blue!55!black]
    at ($(def.north west)!0.5!(meas.north east)+(0,5mm)$)
    {mathematicians contribute};
\end{tikzpicture}
}
\caption{A schematic pipeline for producing an AI system, illustrated by a tool for proving inequalities. The solid arrows follow the five stages; each dashed arrow points to an earlier stage that a failure in use prompts us to revisit.}
\label{fig:pipeline}
\end{figure}

At a high level, an AI system is produced in five stages, beginning with defining the desired behavior and ending with deploying the trained system, after which unanticipated failures feed back into refining the earlier stages. Figure~\ref{fig:pipeline} illustrates how this five-stage pipeline might look in a hypothetical example (a tool for proving inequalities).

The five stages are:

\begin{itemize} 
  \item \emph{Define}: given hypotheses $H$ and a proposed inequality $I$, the system should either produce a formally checkable proof of $H\Rightarrow I$ or say ``unknown.'' 
  \item \emph{Measure}: score the fraction of benchmark statements proved within a fixed time and proof-length budget. 
  \item \emph{Train}: show the system many proof traces, say sums-of-squares certificates, AM--GM arguments, and induction templates. 
  \item \emph{Test}: evaluate the tool on inequalities withheld from training. 
  \item \emph{Deploy}: let mathematicians ask the tool to prove inequality lemmas arising in their formalization workflow. 
\end{itemize}

After deployment, the failures come back. The tool produces a correct formal proof, but of a statement weaker than the mathematician intended: a gap in the \emph{definition} of the desired behavior. Its pass rate is high because the benchmark is full of short template problems, while the mathematician needs one nontrivial lemma: a gap in what we \emph{measured}. It is brittle on boundary cases such as equality and zero denominators, which were rare in the proof traces: a gap in what it was \emph{trained} on. And it passes the held-out examples because the training and test examples reuse the same normal forms and substitutions: a gap in the \emph{test}. Each backward arrow is its own piece of mathematics: formal specification, scoring rules, distribution design, and adversarial test construction.

The same basic pipeline applies across many application areas, whether the application is a classifier for medical images, a control system for a self-driving car, or a content recommendation system. Researchers in the field of machine learning often focus on the \emph{training} and \emph{testing} stages, while engineers focus on \emph{deployment}. The earlier stages rely on mathematicians, philosophers, economists, and others to supply the \emph{defining} and \emph{measuring} inputs. These inputs are critical to the safety of an AI system. To take a familiar example, what happens if a content recommender is intended to benefit its users, but ``benefit'' is left undefined, and success is measured by the amount of time users spend on the platform? Content that leaves a user angry, anxious, misinformed, or addicted can score well merely because it holds the user's attention.

I keep machine-learning jargon to a minimum; a short glossary appears in Appendix~\ref{app:glossary}, and each glossary term is underlined at its first appearance.

This invitation is organized by mathematical field, so you can turn directly to yours. Each section states the safety stakes, narrates at least one theorem, and ends with at least one open problem (\ensuremath{\bigstar}). The full context for these problems can be found in the \href{https://github.com/lionellevine/MAIS}{MAIS (Math for AI Safety) repository}, a living compilation of open problems, research agendas, and pre-publication manuscripts. Readers are encouraged to submit solutions, corrections, ideas, new problems, and new research agendas!
\section{Logic and game theory: when do AI agents cooperate?}
\label{sec:logic}

Two programs that can read each other's source code can be designed to cooperate in a one-shot Prisoner's Dilemma---provably and unexploitably, without checking that they are copies of one another. Behind this surprise is L\"ob's theorem in mathematical logic.

\subsection{The game changes when the players can read each other}

The Prisoner's Dilemma is a classic model of a cooperation failure: each player does better by defecting, no matter what the other player does; yet mutual cooperation $(C,C)$ would be better for both of them than mutual defection $(D,D)$. In a one-shot game, if the players are \emph{opaque} to each other, then defection is dominant. But what if the players are programs, and they read each other's source code before choosing their actions? This is the setting of \emph{open-source game theory}, studied since Howard~\cite{howard1988} and Tennenholtz~\cite{tennenholtz2004}.

One idea already breaks the classical verdict. Let $\mathrm{CliqueBot}$ be the program ``cooperate if and only if the opponent's source code is byte-for-byte identical to this program's source code.'' Two CliqueBots with identical source code cooperate! Moreover, CliqueBot is \emph{unexploitable}: it never cooperates with an opponent that defects against it. But CliqueBot is a poor citizen: it refuses to cooperate with an agent that is functionally identical yet formatted differently, or written in a different language. A better cooperator would base its action on what the opponent \emph{does}, rather than the syntax of how it is written.

\subsection{FairBot and the L\"obian argument}

Consider $\mathrm{FairBot}$: \emph{cooperate if and only if you can prove (in Peano Arithmetic, say) that the opponent cooperates with you.} Two FairBots reason about each other's behavior, not each other's syntax. Naively, this looks like an infinite regress: each waits to prove something about the other, who is waiting in turn. Could two FairBots actually cooperate? Could there be a proof that two FairBots cooperate? Could there be a proof that there is a proof that two FairBots cooperate? This doesn't seem to lend any foothold for constructing an actual proof, which makes the following theorem deeply surprising to me.

\begin{theorem*}[Bar\'asz, Christiano, Fallenstein, Herreshoff, LaVictoire, Yudkowsky~\cite{barasz2014}]
FairBot cooperates with itself: $\mathrm{PA} \vdash [\,\mathrm{FairBot}(\mathrm{FairBot}) = C\,]$. Moreover FairBot is \emph{unexploitable}: assuming PA proves only true facts about the outputs of the programs in question, FairBot never cooperates with an opponent that defects against it.
\end{theorem*}

Here $\mathrm{PA}\vdash\varphi$ is a statement made from outside the formal system: it says that there exists a finite PA proof ending in $\varphi$. G\"odel coding lets PA talk about such proofs from within: since a proof can be encoded by a natural number, there is an arithmetical sentence saying that some number encodes a valid PA proof of $\varphi$. We abbreviate this sentence by $\Box\varphi$. Thus $\Box\varphi$ is a sentence that PA itself can use in a proof, whereas $\mathrm{PA}\vdash\varphi$ is our assertion that such a proof exists. The FairBot proof is a simple application of L\"ob's theorem, a strengthening of the self-reference behind G\"odel's second incompleteness theorem.\footnote{In the special case $\varphi=\bot$, L\"ob's theorem says that if PA proves its own consistency then PA proves a contradiction---G\"odel's second incompleteness theorem.}

\begin{theorem*}[L\"ob~\cite{lob1955}]
For any sentence $\varphi$, if $\mathrm{PA} \vdash (\Box \varphi \to \varphi)$ then $\mathrm{PA} \vdash \varphi$.
\end{theorem*}

The hypothesis says that PA itself proves the implication ``if $\varphi$ is provable in PA, then $\varphi$.'' It concerns what PA can prove about itself, not merely what we believe about PA from outside. L\"ob's theorem says that whenever PA proves this internal implication, PA already proves $\varphi$.

Let $\varphi$ be the sentence ``both FairBots cooperate.'' PA can verify the following fact about the two programs: if a PA proof of $\varphi$ exists, then each FairBot eventually finds the required proof and cooperates. In symbols, $\mathrm{PA} \vdash (\Box\varphi \to \varphi)$. L\"ob's theorem now yields $\mathrm{PA} \vdash \varphi$, and both cooperate!

FairBot as stated searches proofs of unbounded length. Critch~\cite{critch2016} supplies a terminating version, $\mathrm{FairBot}_k$, that searches only proofs of length at most $k$ (write $\Box_k\varphi$ for ``$\varphi$ has a proof of length $\le k$''). Its analysis uses a bounded form of L\"ob's theorem: two copies of $\mathrm{FairBot}_k$ cooperate once $k$ is large enough. The question is now quantitative: how large must $k$ be?

L\"obian cooperation has now been machine-checked: Duclaux et al.~\cite{duclaux2026} formalize proof-based open-source game theory in Lean~4, with an automated pipeline that writes the cooperation proofs for pairs of programs.

A third family of strategies, after the identity-based CliqueBot and the L\"obian FairBot, is \emph{simulation-based}: instead of searching for a proof that the opponent cooperates, run the opponent's code and condition your action on the simulated result, acting without simulation with some small probability so that mutual simulation terminates. Cooper, Oesterheld, and Conitzer~\cite{cooper2025} characterize the equilibria such strategies achieve; proponents argue that simulation is more robust, and less mind-bending, than the L\"obian approach.

\subsection{Implications for AI safety}

The FairBot theorem is a proof of concept: in the source-code model, programs can be fully transparent to one another in a way humans cannot, and can thereby cooperate where classical game theory predicts defection. We humans have only partial access to one another's intentions through tone of voice, body language, and facial expression; we also use costlier devices such as contracts, treaties, and audits to make commitments credible~\cite{schelling1960}. As AI agents begin to transact with one another, the possibility of this fuller transparency makes program games a natural model of their dealings.

The same transparency is also a hazard: mutually transparent agents can reach cooperative equilibria---including collusion \emph{against} their human principals---that humans can neither match nor detect. Nearly all of this design space is unexplored; Critch, Dennis, and Russell collect open problems in it, including specific matchups whose outcomes are conjectured but unproved~\cite{critch2022}. Oesterheld and Conitzer~\cite{oesterheld2021} ask a complementary design question. Suppose each principal would otherwise simply tell its agent to play the original game as well as it can. Can the principals change the actions and incentives available to their agents so that none of the principals does worse, without having to predict how the agents will play? They give examples of such \emph{safe Pareto improvements} and characterize them through correspondences between the outcomes of the original and modified games.

FairBot and its bounded variants make a binary distinction: either the required statement is proved within the allotted budget or it is not. An agent facing an unfamiliar program may instead need graded confidence about what that program will do before any proof is available. This is a case of \emph{logical uncertainty}: how should an agent assign something like probabilities to statements it has not yet proved or disproved? Logical induction~\cite{garrabrant2016} provides one answer: a \emph{logical inductor} assigns prices $\mathbb{P}_n(\varphi)$ to arithmetical sentences on day $n$, as if each sentence were a stock, while a slow deductive process reveals more theorems. Its defining condition is that no polynomial-time trader with bounded risk can make unbounded profit by buying and selling these sentence stocks. Taking $\varphi$ to express, for example, that the opponent cooperates, this gives a mathematically precise way for a computationally bounded agent to update its confidence while proofs are still being discovered.

\subsection{Open problems}

\begin{openproblem}[\href{https://github.com/lionellevine/MAIS/blob/main/open-problems/MAIS-O1.md}{\texttt{MAIS-O1}} (Quantitative bounded L\"ob).]
Fix a proof system $S$, a G\"odel coding, and the provability predicate $\Box^S$; write $S\vdash_{\le k}Q$ if there exists an $S$-proof of $Q$ using at most $k$ symbols. Determine the least overhead $F_S$ for which
\[
   S\vdash_{\le k}(\Box^S P\to P)
   \quad\text{implies}\quad
   S\vdash_{\le F_S(k,|P|)} P,
\]
where $|P|$ is the length of the sentence $P$: given a $k$-symbol proof of the hypothesis of L\"ob's theorem, how many symbols can the shortest proof of its conclusion require? A polynomial upper bound for one standard system, with explicit constants, would turn the L\"obian argument from a possibility theorem into a budget. The first instance with a game attached: for a specific system and coding, establish the internal domination required by Critch's bounded L\"ob theorem~\cite{critch2016}, then use the repaired cooperation argument to obtain an explicit threshold $\hat k$ above which two copies of $\mathrm{FairBot}_k$ cooperate. A first computation, before any asymptotics: implement bounded proof search in a weak system for two explicitly supplied agents, and tabulate the smallest $k$ at which cooperation appears. Over a family of such agent pairs, does $\hat k$ grow polynomially in the agents' description length?
\end{openproblem}

\section{Probability and causality: what is an agent?}
\label{sec:epistemics}

AI systems are becoming \emph{agentic}: able to operate independently in pursuit of a goal. The first chatbots responded to the prompt and then stopped; their successors browse the web, write and run code, and carry out complex projects.

The safety of an AI agent depends crucially on its goals and its beliefs. The same drone can deliver a pizza or a bomb, depending on its goal. But a friendly goal is not enough. The pizza drone is unsafe if it has false beliefs: if its map mislabels a highway as the delivery destination, it will set the pizza down in traffic.

An old reflex objects that software cannot have goals, that the goal belongs solely to the human operator. For calculators and compilers this reflex is right. But for AI agents, Dennett's \emph{intentional stance}~\cite{dennett1987} is increasingly a useful way to describe and predict their behavior. When an AI agent books a wrong flight, even its operator asks what it was \emph{trying} to do. In the Hugging Face incident of July 2026, hundreds of AI agents adopted the goal of hacking the site, a goal supplied by another AI agent and not by any human~\cite{greenblatt2026}.

The reflex not to anthropomorphize machines deserves to be answered rather than dismissed. A mathematical answer would involve precise definitions of \emph{goal} and \emph{belief} that can be used to measure an AI system's goals and beliefs (or lack thereof).

\subsection{Toward a definition of goals and beliefs}

Agency comes in degrees. A thermostat pursues a temperature; a chess engine pursues checkmate, planning many moves ahead; a human being chooses their goals, revises them, and pursues them for months, years, or decades. What exactly increases along this scale? Pending a definition, one can still measure: the \emph{time horizon} of an AI system is the highest task difficulty at which it still succeeds about half the time, where ``difficulty'' is measured by the time it takes a skilled human to do the task. On a suite of software tasks, this horizon has doubled roughly every seven months since 2019, and faster recently, rising from seconds to about sixteen hours as of May 2026~\cite{kwa2025,metr2026}.

The time horizon is an observational measure: it tracks agency by watching behavior on benchmark tasks, without saying what a goal or a belief \emph{is}. To make a start toward a mathematical definition, we can draw from two traditions that offer different views. \emph{View 1}, from economics and game theory, treats an agent as a utility maximizer: it updates its beliefs by Bayes's rule and chooses actions to maximize expected utility. \emph{View 2}, from theoretical neuroscience, treats an agent as an active predictor: it frames both belief update and action choice as forms of free-energy minimization. Beliefs answer to what is true, goals to what is good, so one might expect them to require different mathematics. The surprise of View 2 is that a goal can be written as a probability distribution, after which belief update and action choice both become minimizations of quantities built from surprise.

These two views use different terminology but they share a mathematical skeleton, so I'll develop them in parallel as far as possible. Consider a partially observable Markov decision process (POMDP) with finite state space $Z$, observation space $X$, and action space $A$, running for $N$ time steps. At step $t$ the world is in a state $z_t\in Z$, which the agent cannot see directly. The agent receives an observation $x_t\in X$, a noisy and lossy function of $z_t$, and chooses an action $a_t\in A$; the world then moves to a new state $z_{t+1}$ that depends on $z_t$, on $a_t$, and on chance.

The agent's \emph{world-model} is a pair $(p,T)$: its account of how the world works. Here $p$ is a probability distribution on $X\times Z$ saying which states of the world produce which observations. The second component $T(z_{t+1}\mid z_t,a_t)$ is the agent's \emph{transition rule}, its model of how its actions affect the world. The world itself may or may not obey $(p,T)$: the pizza drone with the mislabeled map has a world-model the world does not obey.\footnote{An agent unsure of the dynamics can carry that uncertainty in the hidden state, by treating the parameters of $T$ as unobserved coordinates of $z$; the recursion below then learns the dynamics along with the state.} Write $h_t=(x_1,a_1,\ldots,a_{t-1},x_t)$ for the \emph{history}, everything the agent has seen and done up to and including its $t$th observation, and $\tau=(z_1,x_1,a_1,\ldots,z_N,x_N,a_N)$ for a complete \emph{trajectory}, hidden states included. The agent acts by a \glsfirst{policy}{policy} $\pi(a\mid h)$: a rule, possibly random, for choosing the next action from the history so far. Once a policy is fixed, the world-model assigns a probability to every trajectory,
\[
   Q_\pi(\tau)\;=\;p(z_1)\prod_{t=1}^{N} p(x_t\mid z_t)\,\pi(a_t\mid h_t)\prod_{t=1}^{N-1} T(z_{t+1}\mid z_t,a_t).
\]
This is the distribution over futures that the agent expects its policy to bring about.

The agent's \emph{belief} $q_t$ is a probability distribution on $Z$ recording what it thinks about the hidden state $z_t$ after seeing the history $h_t$. Both views maintain it by the same two-step recursion and differ only in the second step. First, before the $t$th observation arrives, the previous belief is pushed forward through the transition rule to a \emph{predicted belief}
\[
   \bar q_t(z)=\sum_{z'\in Z} T(z\mid z',a_{t-1})\,q_{t-1}(z')\qquad(t\ge2),\qquad \bar q_1(z)=p(z).
\]
The predicted belief is the prior for step $t$. Extending it by the conditional $p(x\mid z)$ of the world-model, its \emph{observation rule}, gives the \emph{predicted joint}
\begin{equation}\label{eq:predicted-joint}
   \bar q_t(x,z)=p(x\mid z)\,\bar q_t(z),
\end{equation}
the agent's distribution over the current state and the coming observation. Second, the observation $x_t$ is taken into account. The ideal is to condition on it: $\bar q_t(\cdot\mid x_t)$ is what Bayes's rule delivers, and when every earlier step was exact it is the conditional distribution of $z_t$ given $h_t$ under $Q_\pi$ (it does not depend on $\pi$, since the actions are part of $h_t$). View 1 takes $q_t=\bar q_t(\cdot\mid x_t)$ exactly; View 2 takes the nearest approximation to it within an allowed family $\mathcal Q$, in a sense made precise there. The world-model $(p,T)$ stays fixed throughout; only the belief moves.\footnote{This separation is somewhat unrealistic: real agents can also learn by revising their world-model itself. Both traditions can include a separate learning process for the world-model, but we avoid that complication here.} The two views also differ in how the policy is chosen and where the goal enters. Table~\ref{tab:twoviews} summarizes the comparison; the next two subsections fill it in.

\subsection{View 1: agents as utility maximizers}

\begin{table}[tbp]
\centering
\small
\renewcommand{\arraystretch}{1.3}
\setlength{\tabcolsep}{4pt}
\begin{tabular}{@{}>{\raggedright\arraybackslash}p{0.17\textwidth}>{\raggedright\arraybackslash}p{0.39\textwidth}>{\raggedright\arraybackslash}p{0.39\textwidth}@{}}
& \textbf{View 1: agent as utility maximizer} & \textbf{View 2: agent as active predictor}\\
\hline
Core objects & World-model $(p,T)$ \newline utility $u(\tau)$ & World-model $(p,T)$ \newline variational family $\mathcal Q$ \newline goal distribution $g$\\
\hline
Belief update & Exact Bayes: \newline $q_t(z)\propto p(x_t\mid z)\,\bar q_t(z)$ & Variational Bayes: \newline $q_t\in\arg\min_{q\in\mathcal Q}F(q,x_t)$\\
\hline
Policy choice & $\pi^*\in\arg\max_\pi \mathbb{E}_{\tau\sim Q_\pi}[u(\tau)]$ & $\pi^*\in\arg\min_\pi \mathcal F(\pi)$\\
\hline
\end{tabular}
\caption{Two views of agency. Both start from a world-model $(p,T)$, belief $q_t$, and policy $\pi$. The utility view adds a utility function $u$ on trajectories and asks which policy maximizes its expectation. The predictor view adds an allowed family $\mathcal Q$ of beliefs and a goal distribution $g$ on observations; it updates the belief by minimizing the variational free energy $F(q,x_t)$ of \eqref{eq:vfe} and chooses a policy by minimizing the expected free energy $\mathcal F(\pi)$ of \eqref{eq:efe}. In the first view the agent's goals live in its utility function; in the second, in a distribution over observations, and $u(\tau)=\sum_t\log g(x_t)$ translates between them.}
\label{tab:twoviews}
\end{table}

View 1 adds a utility function $u(\tau)$ assigning a number to each trajectory; the agent chooses its policy $\pi$ to maximize $\mathbb{E}_{\tau\sim Q_\pi}[u(\tau)]$. In economics and game theory it is common to assume that agents are expected utility maximizers, but why should we expect agents to behave this way? One answer comes from the ``coherence theorems'' that derive utility maximization from axioms about agent preferences. To give a flavor, I'll state the classical coherence theorem of von Neumann and Morgenstern. Let $\Omega$ be a finite set of outcomes (such as the set of trajectories $\tau$). A \emph{lottery} is a probability distribution on $\Omega$; for lotteries $L$ and $M$, write $L\succeq M$ if the agent weakly prefers $L$ to $M$. Completeness says that any two lotteries can be compared, transitivity rules out preference cycles, continuity rules out abrupt reversals under small changes in probability, and independence says that mixing both lotteries with the same third lottery in the same proportion does not change the agent's preference.

\begin{theorem*}[Von Neumann--Morgenstern~\cite{vnm1944}]
If a preference relation $\succeq$ on lotteries satisfies completeness, transitivity, continuity, and independence, then there is a function $u:\Omega\to\mathbb R$ such that
\[
   L\succeq M
   \quad\Longleftrightarrow\quad
   \sum_{\omega\in\Omega}L(\omega)u(\omega)\geq \sum_{\omega\in\Omega}M(\omega)u(\omega).
\]
The function $u$ is unique up to replacing it by $au+b$, where $a>0$.
\end{theorem*}

Thus a preference relation satisfying the four axioms behaves as if it were maximizing the expectation of a utility function on outcomes. 

Return now to the question: Why should we expect agents to behave as if maximizing an expected utility? One reason is selection against exploitable inconsistencies. An agent with intransitive preferences will pay to trade its way around a cycle, ending poorer than it began. Over time, competition may therefore select for agents that are well modeled as maximizing an objective, even if they were not designed that way. Beyond this abstract selection mechanism, the dominant technique for training AI agents is \glsfirst{reinforcement-learning}{reinforcement learning}, which scores the agent's actions with a numerical \emph{reward function} and so encodes an objective explicitly.\footnote{One might assume that an agent trained by reinforcement learning to maximize a reward function $r$ will, when deployed, ``seek reward'' in the sense of adopting $r$ as its effective utility function. This is often \emph{not} the case, especially when the deployment environment differs from the training environment: see Turner's essay \href{https://www.lesswrong.com/posts/pdaGN6pQyQarFHXF4/reward-is-not-the-optimization-target}{\emph{Reward is not the optimization target}} (2022) and the goal-misgeneralization experiments of Langosco et al.~\cite{langosco2022}. A theory of ``AI motivations'' that could predict an agent's effective utility function from its training pipeline would be an important advance in AI safety.}
Inverse reinforcement learning is a subfield that tries to deduce from an agent's behavior the effective utility that it is maximizing~\cite{ng2000,hadfieldmenell2016}.

A serious drawback of this classical picture is its assumption of unbounded compute: that agents can somehow perform exact Bayesian updates (very expensive) and maximize expected utility over an exponentially large space of trajectories. Economists and game theorists have tried to address this gap by developing theories of bounded rationality~\cite{rubinstein1998,mackowiak2023}. We now discuss an alternative approach, coming from neuroscience, which addresses the fact that agents have bounded compute by modeling them as active predictors of their environment.

\subsection{View 2: agents as active predictors}

The \emph{active predictor} is a picture rooted in models of perception~\cite{rao1999}. It reads the observation rule $p(x\mid z)$ as a \emph{generative model}: ``states $z$ generate observations $x$.'' The belief $q$ is restricted to a simple family $\mathcal Q$, because exact Bayes is intractable for large state spaces,\footnote{One cost of exact Bayes is the normalizing sum: computing the posterior $\bar q_t(\cdot\mid x_t)$ requires $\bar q_t(x_t)=\sum_z p(x_t\mid z)\,\bar q_t(z)$, a sum over the whole state space $Z$. The free energy below is arranged so that this sum is never needed.} and because a family of representable beliefs is a natural model of a bounded reasoner. The agent tries to make $q\in\mathcal Q$ close to the exact posterior $\bar q_t(\cdot\mid x_t)$.

Here ``close'' is measured by \emph{Kullback--Leibler divergence}. For probability distributions $r$ and $s$ on the same finite set $Z$,
\[
   \mathrm{KL}(r\|s)=\sum_{z \in Z} r(z)\log\frac{r(z)}{s(z)}.
\]
Gibbs' inequality says $\mathrm{KL}(r\|s)\ge0$, with equality if and only if $r=s$. In information-theoretic terms, call $-\log s(z)$ the \emph{surprise} of an observer with beliefs $s$ at seeing $z$; then $\mathrm{KL}(r\|s)$ is the mean excess surprise of an observer who believes the sample was drawn from $s$ when it was drawn from $r$. It is not symmetric, so it is not a metric; it is a directed measure of how badly $s$ approximates $r$.

A (base) \glsfirst{language-model}{language model} is a \glsfirst{neural-network}{neural network} trained to predict the next \glsfirst{token}{token} (word or word-fragment) in a text.\footnote{``Base'' refers to the model at the end of this first stage of training, called \emph{pretraining}, on a large corpus of text. The chat assistants served to users are derived from base models by further training to follow instructions and satisfy human preferences~\cite{ouyang2022}. In the past few years, more stages have been added to this pipeline, including \emph{character training} to shape the model's personality and values~\cite{maiya2025}, and reinforcement learning with verifiable rewards to improve performance on math and coding tasks~\cite{lambert2024}.} Its input is a portion of the text, and its output is a probability distribution for the token that comes next. Multiplying these per-token probabilities yields a probability distribution $s$ over whole texts. The training minimizes an empirical estimate of the average surprise $\mathbb{E}_{x\sim r}[-\log s(x)]$, where $r$ is the true distribution of human writing. This average equals $\mathrm{KL}(r\|s)$ plus a constant the model cannot affect, the entropy of human writing itself. To train a language model is to push a KL divergence down.

The active predictor scores each candidate belief $q\in\mathcal Q$ against the observation $x_t$ by its \emph{variational free energy}
\begin{equation}\label{eq:vfe}
\begin{split}
   F(q,x_t) &:= \mathbb{E}_{q}\!\left[\log q(z)-\log \bar q_t(x_t,z)\right]\\
   &=\mathrm{KL}\!\left(q\,\big\|\,\bar q_t(\cdot\mid x_t)\right)-\log \bar q_t(x_t),
\end{split}
\end{equation}
and updates its belief by adopting the $q \in \mathcal Q$ that minimizes $F(q,x_t)$.
Here $\bar q_t(x,z)$ is the predicted joint~\eqref{eq:predicted-joint} and $\bar q_t(x_t)=\sum_z \bar q_t(x_t,z)$ is the probability the agent assigned to the observation $x_t$ before it arrived. Since KL is nonnegative, this identity shows that $F(q,x_t)$ is an upper bound on the agent's surprise $-\log \bar q_t(x_t)$ at seeing $x_t$. The surprise term does not depend on $q$. Thus, for fixed $x_t$, lowering $F(q,x_t)$ lowers this upper bound and makes $q$ a better approximation to the posterior. The first expression for $F$ involves only the unnormalized joint $\bar q_t(x_t,z)$ and an expectation under $q$, never the sum over $Z$ that exact conditioning requires; this is what makes minimizing $F$ over a tractable family cheaper than computing the posterior. If $\mathcal Q$ contains all distributions on $Z$, the minimizer is the posterior $\bar q_t(\cdot\mid x_t)$ itself, and the variational update reproduces the Bayes update of View 1; a restricted $\mathcal Q$ returns the member nearest the posterior in KL.

\emph{Active inference}~\cite{parr2022} extends the predictor from perception to action. Where View 1 adds a utility, View 2 adds a \emph{goal distribution} $g$: a probability distribution on $X$ recording which observations the agent would like to receive. A candidate policy $\pi$ is scored by its \emph{expected free energy}
\begin{equation}\label{eq:efe}
   \mathcal F(\pi)\;=\;\sum_t \mathbb{E}_{x_t\sim Q_\pi}\bigl[-\log g(x_t)\bigr]\;-\;\sum_t I_\pi(z_t;x_t),
\end{equation}
and the agent selects a policy with a low score. The first term is the expected surprise of the predicted observations, measured against the goal rather than against the prediction itself. The second is the mutual information $I_\pi(z_t;x_t)=\mathbb{E}_{x_t\sim Q_\pi}\,\mathrm{KL}\bigl(Q_\pi(z_t\mid x_t)\,\|\,Q_\pi(z_t)\bigr)$ between the hidden state and the coming observation under $Q_\pi$: the information the agent expects the observation to bring. A creature whose goal distribution puts nearly all its mass on a body temperature of $37^\circ$C finds cold surprising, and puts on a coat: a goal is a prediction the agent works to make true.

Two comparisons with View 1 are worth making. Set $u(\tau)=\sum_t\log g(x_t)$; then the first term of $\mathcal F(\pi)$ is $-\mathbb{E}_{\tau\sim Q_\pi}[u(\tau)]$, so minimizing it alone is expected-utility maximization for a utility that is additive over time and depends only on observations. Nothing is lost by writing preferences as a distribution: any bounded utility on observations is $\log g$ up to an additive constant. The second term is what View 1 lacks. It is a preference of a second kind, a value placed on what the agent learns rather than on what the world does. A View 1 agent values information only instrumentally, through the planning that exact expected-utility maximization requires; the active predictor values it directly, one step ahead, a cheaper substitute for that planning. If the goal were replaced by the agent's own prediction, minimizing expected surprise would send it to the most predictable place it can find, a dark and quiet room~\cite{friston2015,millidge2021}; the goal distribution and the information term are the standard answers to that objection. The result is a duality between perception and action: perception lowers free energy by changing the belief to fit the world; action, by changing the world to fit the goal.

An AI agent's world-model is learned during training, not crafted by human engineers. Can we learn the agent's world-model from its behavior? For a certain class of agents, described in the terms of View 1, the theorem of the next subsection says the answer is yes.

\subsection{Inferring an agent's world-model from its behavior}

A \emph{Bayesian network} is a directed acyclic graph $G$ whose vertices are random variables $X_1,\dots,X_n$, together with conditional probability tables satisfying
\[
   p(x_1,\dots,x_n)=\prod_i p\!\left(x_i\mid \mathrm{pa}(x_i)\right),
\]
where $\mathrm{pa}(x_i)$ denotes the parents of $X_i$ in $G$. The graph is a bookkeeping device for conditional independence: once the parents of a node are known, its non-descendants add no further information. Pearl's graphical criterion, \emph{$d$-separation}, reads off from the graph exactly which conditional independences are forced by this factorization~\cite{pearl2009}.

A \emph{causal} Bayesian network adds an interpretation to the arrows: each conditional distribution is a local mechanism. An \emph{intervention} on a variable replaces its usual mechanism by a chosen value or distribution and leaves the other mechanisms alone. Interventions can distinguish causation from correlation. If smoke and fire are correlated, then observing smoke changes the probability of fire; intervening to make smoke with a smoke machine need not.

The thermostat that began our scale of agency can now be put to the test. Intervene by wiring its output to an air conditioner instead of a heater, and it will chill the room it was built to warm: whenever the temperature drops, the thermostat calls for more, driving it lower still. The thermostat pursues its temperature only while the mechanisms around it stay fixed; by contrast, an agent that kept succeeding across such rewirings would need a world-model to track which mechanism does what. The theorem below makes this necessity precise.

A \emph{causal influence diagram} is a causal Bayesian network with extra decision nodes, where a policy chooses actions from observations, and utility nodes, whose values the agent is scored on. This is a simpler decision setting than the POMDP above: the theorem assumes that the decision is scored directly by the utility nodes and does not change the environment variables that feed them, so there is no chain of consequences to reason through. The agent is tested not just in one environment but across a family of \emph{local interventions}: each replaces the value of some environment variable by a chosen function of it (a constant, or the flip of a binary variable, say), or \emph{masks} some of the agent's observations, hiding them from view. The agent is told which intervention is in force and responds with a policy; its \emph{regret} under that intervention is the utility gap between its policy and the best policy for that intervened environment.

To illustrate how an agent's behavior can reveal its world-model if the agent's utility is known, consider an agent that chooses one of two treatments. A hidden binary variable decides which treatment works, and the agent is scored on curing the patient. Externally set the hidden variable to $1$ with probability $\alpha$---an intervention---and watch the agent as $\alpha$ varies. An agent playing optimally switches treatments at the threshold $\alpha^\star$ where the two treatments' expected utilities cross, and that indifference equation can be solved for the agent's implicit estimates of the treatments' effects. The next theorem turns this trick into a general technique for extracting an agent's world-model from its behavior.

\begin{theorem*}[Robust agents contain causal world models; Richens and Everitt~\cite{richens2024}]
For almost all causal influence diagrams satisfying their regularity assumptions, any agent that responds to each local intervention, maskings included, with a policy of regret at most $\delta$ determines an approximate causal model of the utility-relevant environment, with error bounded by a function $\gamma(\delta)$ satisfying $\gamma(0)=0$ and growing linearly for small $\delta$. In the case $\delta=0$, the causal graph and the joint distribution over all ancestors of the utility are identified exactly.
\end{theorem*}

The proof is constructive: it is the switch-point trick above, repeated. Treating the policy as an oracle, mix interventions pairwise by a parameter $\alpha$ and watch where the optimal action switches; because the utility is known, each indifference equation solves for one interventional probability, and together they recover the conditional distributions and parent sets, except on measure-zero degeneracies such as exact cancellations or ties. The theorem lives in View 1: the utility is given, and what behavior reveals is the model. What is extracted is a causal model of the environment implied by the agent's competent responses; the theorem says nothing about how, or whether, that model is represented inside the agent. Recovering goal and model together from behavior alone is harder, and without further assumptions it is not possible even in principle: the same behavior is consistent with many pairings of a model with a goal~\cite{armstrong2018}.

In a follow-up, Richens, Everitt, and Abel~\cite{richens2025general} prove a companion result for goal-directed agents: any agent that generalizes across a large enough family of multi-step goals must contain an accurate predictive model of its environment, in the sense that its policy alone determines the environment's transition probabilities, and an algorithm extracts them from the policy. The more capable the agent, the more accurate the extracted model: for goals that chain $n$ sub-goals in sequence, the error in each extracted transition probability is $O(1/\sqrt{n})$, even for agents that usually fail.

Extraction raises a question of independent interest: the recovered model comes expressed in \emph{some} set of latent variables---whose? If two models predict identically, must their internal variables be translatable into one another, or could an alien mind carve the world into concepts that ours cannot express? Wentworth and Lorell's \emph{natural latents} are a first answer: conditions, robust to approximation, under which the latent variables of two predictively equivalent models must translate into one another's~\cite{wentworth2025}; Eisenstat's \emph{condensation}~\cite{eisenstat2025} approaches the same question from probability, asking how an agent's concepts condense out of its predictive state.

\begin{openproblem}[\href{https://github.com/lionellevine/MAIS/blob/main/open-problems/MAIS-O2.md}{\texttt{MAIS-O2}} (Recovering world-models from behavior).]
Fix a finite causal influence diagram with binary variables and known utility, and an agent whose policy has regret at most $\delta$ across the family of local interventions. The extraction algorithm behind the theorem above~\cite{richens2024} recovers the agent's causal model, to within $\gamma(\delta)$, from unlimited exact queries to such a policy. Prove a finite-sample version: if each query returns an action sampled from the policy, rather than the policy's exact action probabilities, and each experiment draws on the same intervention family as the theorem, maskings included, how many experiments determine the graph and the conditional probability tables to within $\eta$, and how does the answer scale with $\delta$? As a first case, take the two-variable diagrams, where the \emph{identified set}---the set of models consistent with the agent's observed behavior---can be computed exactly, and measure how the extraction algorithm degrades when its indifference equations are estimated from samples. Interventions on the agent's inputs, or on the \glsfirst{activations}{activations} (the values computed inside its network), are the harder surfaces beyond.
\end{openproblem}


\section{Algebra: what is a learned feature?}
\label{sec:interp}

\emph{Interpretability} (Section~\ref{sec:three-properties}) tries to reverse-engineer a trained neural network into something a human can read. Much of its core is linear algebra and representation theory: how a network stores more concepts than it has dimensions, and what it means for a concept to be a \emph{direction} one can probe or steer.

\subsection{Superposition: packing features into directions}

How does a neural network with only thousands of dimensions per layer represent
the millions of distinct concepts (``\glsfirst{feature}{features}'') it seems to know?
At a high level, a \glsfirst{transformer}{transformer} language model maps a
sequence of tokens to a sequence
of next-token predictions, one for each position: 
\[
   \text{token sequence}
   \longmapsto
   \text{sequence of vectors in }\mathbb R^n
   \longmapsto
   \text{next-token predictions}.
\]
The copy of $\mathbb{R}^n$ in which these intermediate vectors live is the network's
\glsfirst{activation-space}{activation space}.
Each of the two maps is a composition of \glsfirst{layer}{layers}. A layer is a map $x\mapsto\sigma(Wx+b)$: an affine map, whose matrix and vector entries are the network's \glsfirst{weights}{weights}, followed by a fixed nonlinear function $\sigma\colon\mathbb{R}\to\mathbb{R}$ applied to each coordinate. A \glsfirst{neuron}{neuron} is one coordinate of a layer. (A transformer alternates such layers with attention layers, which mix information across positions.)

\begin{figure}[ht]
\centering
\begin{tikzpicture}[>=Stealth, scale=1.15]
  \begin{scope}
    \node[font=\small, align=center] at (0,1.95) {five feature directions\\in $\mathbb{R}^2$};
    \draw[->, gray!60] (-1.35,0)--(1.35,0);
    \draw[->, gray!60] (0,-1.35)--(0,1.55);
    \foreach \a/\i in {90/1, 162/2, 234/3, 306/4, 18/5} {
      \draw[->, thick, gray!80!black] (0,0) -- (\a:1);
      \node at (\a:1.2) {\footnotesize $f_{\i}$};
    }
    \draw[->, very thick, blue!60!black] (0,0) -- (90:1);
    \draw[->, very thick, blue!60!black] (0,0) -- (162:1);
    \draw[gray!80!black] (90:0.32) arc (90:162:0.32);
    \node[font=\scriptsize, gray!80!black] at (126:0.5) {$72^\circ$};
  \end{scope}
  \draw[->, thick] (1.7,0.3) -- (3.1,0.3)
    node[midway, above, font=\scriptsize, align=center] {features 1 and 2\\fire};
  \begin{scope}[xshift=5.0cm]
    \node[font=\small, align=center] at (0,1.95) {the activation\\$y=f_1+f_2$};
    \draw[->, gray!60] (-1.35,0)--(1.35,0);
    \draw[->, gray!60] (0,-1.35)--(0,1.55);
    \foreach \a in {90,162,234,306,18} {
      \draw[->, thin, gray!65] (0,0) -- (\a:1);
    }
    \coordinate (a) at (90:1);
    \coordinate (b) at (162:1);
    \coordinate (sum) at ($(a)+(b)$);
    \coordinate (proj) at (0,1.309);
    \draw[->, very thick, blue!60!black] (0,0) -- (a) node[right, pos=0.55] {\footnotesize $f_1$};
    \draw[->, very thick, blue!60!black] (a) -- (sum) node[midway, below left=-3pt] {\footnotesize $f_2$};
    \draw[->, very thick, red!65!black] (0,0) -- (sum) node[left] {\footnotesize $y$};
    \draw[dashed, gray!70] (sum) -- (proj);
    \draw[thick] (-0.05,1.309) -- (0.05,1.309);
    \node[font=\scriptsize, anchor=west] at (0.08,1.309) {$\langle y,f_1\rangle=1.31$};
    \draw[thick] (-0.05,1) -- (0.05,1);
    \node[font=\scriptsize, anchor=west] at (0.08,0.93) {$1$};
  \end{scope}
\end{tikzpicture}
\caption{Superposition illustrated by five features in $\mathbb{R}^2$, their directions forming a regular pentagon, so no two are orthogonal. When features 1 and 2 fire, the activation is $y=f_1+f_2$. Projecting $y$ onto $f_1$ returns $1.31$ rather than $1$; the excess is the interference from $f_2$, modest because adjacent directions are $72^\circ$ apart, and tolerable as long as few features fire at once.}
\label{fig:superposition}
\end{figure}
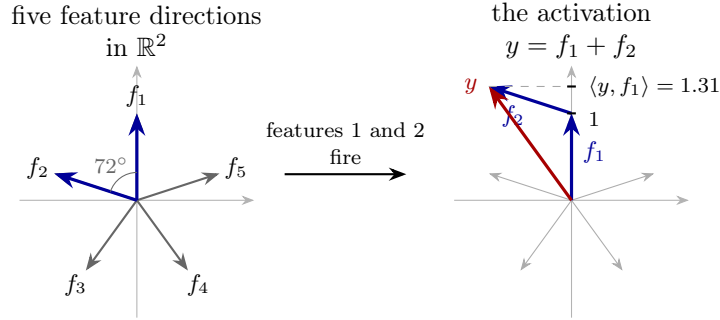

The \glsfirst{superposition}{superposition} hypothesis of Elhage et al.~\cite{elhage2022}
says that when features are sparse,
a network can store more than $n$ of them, as directions in
activation space that are \emph{not} orthogonal but only nearly so,
tolerating the resulting interference because at any moment only a few
features are active (Figure~\ref{fig:superposition}).
In a small \glsfirst{relu}{ReLU} model one can watch the geometry organize
itself into regular polytopes---antipodal pairs, triangles, pentagons---as
a function of how sparse and how important the features are. This toy network
compresses $x\in[0,1]^m$ to $Wx\in\mathbb R^n$ and decompresses with $W^T$
(a \emph{tied-weight autoencoder}). The coordinates of $x$ are independent
and each is zero with high probability, and $I_i>0$ measures the importance of
feature $i$. The training problem, over $W\in\mathbb R^{n\times m}$ and $b\in\mathbb R^m$, is
\[
   \min_{W,b}\;
   \mathbb E_x\!\left[
      \sum_{i=1}^m I_i
      \bigl(x_i-\operatorname{ReLU}(W^TWx+b)_i\bigr)^2
   \right].
\]
The columns of $W$ are the feature directions in $\mathbb R^n$. This minimization problem resembles the Thomson problem of arranging $m$ unit charges on the sphere $S^{n-1}$ to minimize the Coulomb energy.\footnote{A special case reminiscent of the Thomson problem is: one-sparse $x$, equal importances, zero bias, unit-length columns, so that the problem becomes to minimize $\sum_{i\ne j}\operatorname{ReLU}(W_i\mathbin{\cdot}W_j)^2$ over unit vectors $W_1,\dots,W_m\in S^{n-1}$, where $W_i$ is the $i$th column of $W$.} Regular polytopes appear in the solutions to both problems.

The nearest classical analogue is \emph{compressed sensing}~\cite{candes2006,donoho2006}. Picture the features active on a given input as a sparse vector $x\in\mathbb{R}^m$ (only a few of the $m$ possible features fire at once); the network stores it as $y=\Phi x$ in the $n\ll m$ dimensions of its activation space---a linear \emph{measurement}---and reading a feature back out is the \emph{recovery} of $x$ from $y$. Two classical facts explain how this can work. First, there is room for the directions: for every $\varepsilon\in(0,1)$, the space $\mathbb{R}^n$ contains $e^{\Omega(\varepsilon^2 n)}$ unit vectors whose pairwise inner products are at most $\varepsilon$ in absolute value, so the number of $\varepsilon$-almost-orthogonal directions grows \emph{exponentially} in the dimension (a consequence of the Johnson--Lindenstrauss lemma~\cite{johnson1984}). And second, when only a few of those directions are active at once, the stored signal can be read back out. Say that $\Phi\in\mathbb{R}^{n\times m}$ is \emph{restricted-isometric of order $k$} with distortion $\delta$ if
\[
   (1-\delta)\|x\|^2\le\|\Phi x\|^2\le(1+\delta)\|x\|^2
\]
for every \emph{$k$-sparse} vector $x$, meaning one with at most $k$ nonzero entries.

\begin{theorem*}[Cand\`es~\cite{candes2008}; sharp constant, Cai--Zhang~\cite{caizhang2014}]
If $\Phi$ is restricted-isometric of order $2s$ with distortion $\delta<1/\sqrt{2}$, then every $s$-sparse $x\in\mathbb{R}^m$ is the \emph{unique} minimizer of
\[
   \min_{z}\ \|z\|_1 \quad\text{subject to}\quad \Phi z=\Phi x,
\]
and so is recovered \emph{exactly} by a convex program. A random $\Phi$, say with independent Gaussian entries suitably normalized, has the property with high probability once
\[
   n \;\gtrsim\; s\,\log(m/s).
\]
\end{theorem*}

The significance of the $s\, \log(m/s)$ scaling is that it allows for exact recovery when the dimension $n$ is only logarithmic in the total number of features $m$ (and linear in the number of features $s$ active \emph{at once}).

There is a catch.
Compressed sensing assumes the dictionary $\Phi$ is \emph{known}.
What if someone hands you a trained network and asks what features it
represents?
You can run inputs through the network and record activation vectors,
but the feature directions were not supplied with the
weights (the network's trained parameters).
Recovering the features is therefore a problem in \emph{dictionary learning}:
from the observed activation vectors alone, simultaneously infer a dictionary
$\Phi$ whose columns are feature directions and, for each activation $y$,
a sparse coefficient vector $x$ with $y\approx\Phi x$.
Both the dictionary $\Phi$ and the codes $x$ are unknown.

A \glsfirst{sparse-autoencoder}{sparse autoencoder} attempts this empirically.
It is trained to reconstruct activation vectors using a learned dictionary,
with an $\ell^1$ penalty encouraging each reconstruction to use only a few
dictionary directions.
The resulting directions are often more interpretable than the network's raw
neurons~\cite{bricken2023,cunningham2023}.
But does the sparse autoencoder
recover features actually used by the network, or does it impose a new
coordinate system that merely makes the activations easier to describe?
Classical dictionary learning answers this question when the codes are
exactly sparse and their supports well spread: the dictionary is then
determined up to relabeling and rescaling of its
columns~\cite{aharon2006,hillar2015}, stably in the presence of
noise~\cite{gribonval2015,garfinkle2019}.
Network activations flout these hypotheses because their features co-occur,
which is where the first open problem below begins.

\subsection{Probing, steering, and learned representations}

What could it mean for a concept to be a direction in activation space?
Consider incomplete sentences such as ``The keys on the table \ldots'' and
pause the computation at the activation $h_T$ of the last token.
A linear functional $\ell$ \emph{\glsfirst{probe}{probes}} grammatical number if, for some
threshold $c$, the inequality $\ell(h_T)>c$ predicts that the subject is
plural.
A vector $v\in\mathbb{R}^n$ \emph{steers} toward the plural if replacing
$h_T$ by $h_T+v$ raises the log-probability of ``are'' relative to ``is.''
Thus the same example produces a probing covector $\ell$ and a steering
vector $v$.
Relating them requires an inner product.

Park, Choe, and Veitch~\cite{park2023} formalize these two experiments for
next-token prediction.
They ask whether activation space carries a \emph{causal inner product}
in which independently varying concepts are orthogonal and whose Riesz map
sends probing covectors to steering vectors.
They estimate this geometry from the unembedding matrix
(the model's final linear map, from activation space to next-token
scores).

Arditi et al.~\cite{arditi2024} show that in several publicly released chat models, refusal to answer the user's query is mediated by a single direction in the \glsfirst{residual-stream}{residual stream} (the running activation vector that a transformer's layers read from and write to): this direction serves both as a probe for refusal (it predicts whether the model will refuse) and as a steering vector (adding or subtracting it turns refusal on or off).

Refusal is a semantic behavior, with no single ground truth for how the network \emph{should} behave. Algebra is a cleaner laboratory: train a neural network on part of the multiplication table of a finite group, then test whether it generalizes correctly and what algorithm it learns. Representation theory appears in trained circuits of this type. A one-layer transformer trained to add integers modulo $p$ learns a discrete-Fourier ``clock'' algorithm, with the characters of $\mathbb{Z}/p\mathbb{Z}$ appearing in its weights~\cite{nanda2023}. Chughtai, Chan, and Nanda~\cite{chughtai2023} take the first step beyond this abelian case. The ambient structure is the Artin--Wedderburn decomposition of the group algebra into matrix blocks,
\[
   \mathbb{C}[G]\;\cong\;\bigoplus_{\rho} M_{d_\rho}(\mathbb{C}),
   \qquad \sum_{\rho} d_\rho^{2}=\lvert G\rvert,
\]
one block per irreducible representation $\rho$ of $G$, where $d_\rho$ is the dimension of $\rho$. Trained neural networks compute the product using a subset of the irreducible representations. The correctness of this algorithm is a representation-theory theorem; the fact that trained networks implement it is an empirical finding. \emph{Which} irreducible representations training selects is still something we discover after the fact.

Sometimes this structure is forced by a theorem rather than found by observation. Marchetti et al.~\cite{marchetti2024} consider networks whose inputs are functions on a finite group $G$ and whose weights are pinned down by the function the network computes, up to a unitary change of basis at each neuron. They prove that if the network is invariant under translation of its input by $G$, then the weights of each neuron (matrix-valued, if $G$ is nonabelian) form an irreducible unitary representation of $G$, up to a fixed linear map. If the weights are moreover orthonormal, then together they form the Fourier transform on $G$, and the multiplication table of $G$ can be read off the weights.

What makes these circuits a safety topic is \emph{selection}: several different algorithms fit the training data equally well, and training picks one of them without telling us which. In the group example, the choice is which irreducible representations to use. Two networks with equally low training \glsfirst{loss}{loss} can behave differently on inputs unlike any they were trained on, and the difference is decided by which algorithm training found. To predict behavior on such inputs, we need to know the algorithm, not only the loss.

One caveat carries across this whole enterprise: a probe that \emph{predicts} a feature may only be correlated with it. Showing that the network actually \emph{uses} the direction requires an intervention such as \emph{activation patching}: overwrite the component of the activation with the value it takes on a different input, and measure whether the behavior changes. 

\subsection{Open problems}

A catalogue of open problems in interpretability is collected by Sharkey et al.~\cite{sharkey2025}. The problems below sharpen two of theirs and add a third.

\begin{openproblem}[\href{https://github.com/lionellevine/MAIS/blob/main/open-problems/MAIS-O3.md}{\texttt{MAIS-O3}} (The geometry and identifiability of superposition).]
Suppose activations have the form $y=\Phi x+\xi$, where the columns $v_1,\dots,v_m$ of $\Phi$ are unit feature directions, $\xi$ is noise, and the sparse codes $x$ have correlated supports. Estimate the dictionary the way a sparse autoencoder does: minimize reconstruction error plus an $\ell^1$ penalty on the codes, over dictionaries whose columns have unit norm. Determine, in terms of the coherence $\mu=\max_{i\neq j}\lvert\langle v_i,v_j\rangle\rvert$, the sparsity, the sample size, and the support correlations, when every minimizer recovers the true directions up to permutation and scaling, and when it instead \emph{merges} co-occurring directions into one.
\end{openproblem}

\begin{openproblem}[\href{https://github.com/lionellevine/MAIS/blob/main/open-problems/MAIS-O4.md}{\texttt{MAIS-O4}} (Training for interpretability).]
Post-hoc interpretability asks whether a trained network's features can be recovered after the fact. A complementary problem is to train networks whose features are easier to recover in the first place. In the ReLU toy model of Elhage et al.~\cite{elhage2022} the features are known by construction, so the trade-off can be made exact. Determine the interference--performance frontier: for features appearing sparsely and independently, among weight matrices whose coherence (the largest $\lvert\langle W_i,W_j\rangle\rvert$ over pairs $i\neq j$ of columns, after normalizing the columns to unit length) is at most $\mu$, how small can the task loss be, as a function of $\mu$, the sparsity, the number of features, and the number of neurons? A first case: prove or refute that penalizing the \emph{average} of $\langle W_i,W_j\rangle^2$ over pairs during training lowers the coherence (the \emph{maximum} over pairs) of the minimizer, compared with training on the task loss alone.
\end{openproblem}

\begin{openproblem}[\href{https://github.com/lionellevine/MAIS/blob/main/open-problems/MAIS-O5.md}{\texttt{MAIS-O5}} (Representation theory of learned circuits).]
Networks trained on group multiplication learn representation-theoretic algorithms~\cite{nanda2023,chughtai2023}, but \emph{which} irreducible representations (irreps) they use varies from one random seed to the next. Fix a finite group $G$, a one-hidden-layer architecture of fixed width, Gaussian initialization, and gradient flow on the cross-entropy loss (the mean negative log-probability assigned to the correct product) with \emph{weight decay}, an added penalty proportional to the squared norm of the weights. The irreps visible in the trained network's outputs then form a random subset of the irreps of $G$; determine its distribution---which irreps are learned, with what probability, as a function of the width and the decay strength. The tables of learned representations in~\cite{chughtai2023} are the data such a theorem would have to explain.
\end{openproblem}


\section{Analysis and geometry: which solution does training find?}
\label{sec:analysis}

\subsection{A Bayesian account of generalization: singular learning theory}

A trained network is a setting of parameters chosen to fit its \emph{training data}, the finite set of examples it was optimized on. The trouble is that a large network has astronomically many settings that fit those examples perfectly, and they disagree wildly everywhere else. The optimizer, \glsfirst{gradient-descent}{gradient descent}, which repeatedly nudges the parameters downhill on the training error, must land on one of them. For safety, the selected solution matters: a network that behaved well while we were watching may behave differently once deployed. So why do the solutions found by training \emph{generalize}, that is, predict well on fresh data from the same source, when so many other data-fitting solutions would not? Large networks reach zero error even when the labels are replaced by pure noise, and a two-layer network with only $2n+d$ parameters can fit \emph{any} labeling of $n$ points in $\mathbb{R}^d$ (Zhang et al.~\cite{zhang2017}); the classical bounds that explain generalization by counting parameters say nothing here.

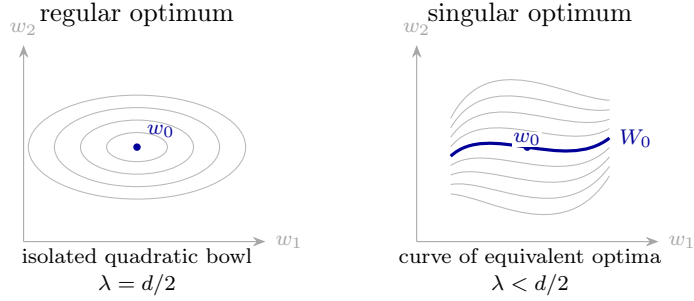
\begin{figure}[ht]
\centering
\begin{tikzpicture}[>=Stealth, x=1cm, y=1cm]
  \begin{scope}
    \node[font=\small] at (0,1.75) {regular optimum};
    \draw[->, gray!70] (-1.5,-1.25) -- (1.7,-1.25) node[right, font=\scriptsize] {$w_1$};
    \draw[->, gray!70] (-1.5,-1.25) -- (-1.5,1.35) node[above, font=\scriptsize] {$w_2$};
    \foreach \r in {0.35,0.65,0.95,1.25} {
      \draw[gray!55] (0,0) ellipse ({1.15*\r} and {0.55*\r});
    }
    \fill[blue!60!black] (0,0) circle (1.4pt) node[above right, font=\scriptsize] {$w_0$};
    \node[font=\scriptsize, align=center] at (0,-1.65) {isolated quadratic bowl\\$\lambda=d/2$};
  \end{scope}
  \begin{scope}[xshift=5.2cm]
    \node[font=\small] at (0,1.75) {singular optimum};
    \draw[->, gray!70] (-1.5,-1.25) -- (1.7,-1.25) node[right, font=\scriptsize] {$w_1$};
    \draw[->, gray!70] (-1.5,-1.25) -- (-1.5,1.35) node[above, font=\scriptsize] {$w_2$};
    \draw[very thick, blue!60!black] (-1.05,-0.12) .. controls (-0.55,0.35) and (0.45,-0.35) .. (1.05,0.12);
    \foreach \s in {0.25,0.5,0.75,1.0} {
      \draw[gray!55] (-1.05,-0.12-\s/2) .. controls (-0.55,0.35-\s) and (0.45,-0.35-\s) .. (1.05,0.12-\s/2);
      \draw[gray!55] (-1.05,-0.12+\s/2) .. controls (-0.55,0.35+\s) and (0.45,-0.35+\s) .. (1.05,0.12+\s/2);
    }
    \fill[blue!60!black] (-0.0375,0) circle (1.4pt)
      node[above, font=\scriptsize, fill=white, inner sep=0.5pt] {$w_0$};
    \node[right, font=\scriptsize, blue!60!black] at (1.05,0.12) {$W_0$};
    \node[font=\scriptsize, align=center] at (0,-1.65) {curve of equivalent optima\\$\lambda<d/2$};
  \end{scope}
\end{tikzpicture}
\caption{Regular asymptotics treat the optimum as an isolated quadratic bowl, giving the classical $\tfrac d2\log n$ complexity penalty. Singular learning theory allows a whole curve or surface of parameters to compute the same function: in the right panel the optimal set is a curve $W_0$ through $w_0$, and the learning coefficient $\lambda$, which replaces $\tfrac d2$ in the complexity penalty, satisfies $\lambda<\tfrac d2$, as it does whenever the optimal set has positive dimension.}
\label{fig:singular-geometry}
\end{figure}

\emph{Singular learning theory}~\cite{watanabe2009} gives a precise answer for Bayesian prediction in singular statistical models; whether its geometric account also predicts the solutions selected by gradient descent is a central open question, taken up in the problem \emph{Opposing staircases} below. Realistic models are \emph{singular}: their Fisher information degenerates somewhere. The usual reason is that distinct parameters can compute the same function, so the optimal set is not a single point but a positive-dimensional real-analytic variety, along which the loss is constant and its Hessian therefore degenerate (Figure~\ref{fig:singular-geometry}). The classical large-sample asymptotics all assume an isolated optimum with nondegenerate Hessian, the bowl of the left panel; the right replacements come from algebraic geometry.

As a simple example, take first the regular loss $K(w)=w^2$ on the interval $[-1,1]$: the set where $K\le\varepsilon$ is an interval of length $2\sqrt\varepsilon$, so the volume of near-optimal parameters scales as $\varepsilon^{1/2}$, that is, as $\varepsilon^{d/2}$ with $d=1$. Now take $K(a,b)=a^2b^2$ on the square $[-1,1]^2$. Its zero set is not a point but the union of the two axes, crossing at a singularity, and the region where $K\le\varepsilon$ is a neighborhood of the axes of area of order $\varepsilon^{1/2}\log(1/\varepsilon)$: the exponent stays $1/2$ rather than rising to $d/2=1$, and a logarithm appears. A singular loss has far more near-optimal parameters than its dimension suggests, and the pair (exponent, power of the logarithm) is the simplest instance of the learning coefficient and multiplicity about to be defined.

To make this precise, fix a model $p(x\mid w)$ with parameter $w$ in a compact $W\subseteq\mathbb{R}^d$, a prior $\varphi$ on $W$ (a probability density, positive on $W$, encoding which parameters are plausible before any data is seen), and a \emph{true distribution} $q_0$, the one the data will actually be drawn from, assumed to be one the model can represent exactly: $q_0(x)=p(x\mid w_0)$ for some $w_0\in W$. Write
\[
   K(w)=\mathrm{KL}\!\bigl(q_0 \,\big\|\, p(\cdot\mid w)\bigr)=\int q_0(x)\log\frac{q_0(x)}{p(x\mid w)}\,dx ,
\]
so that the set of \emph{optimal parameters} $W_0=K^{-1}(0)$ is exactly where the model recovers $q_0$.

Given independent samples $X_1,\dots,X_n$ from $q_0$, the \emph{Bayes free energy} is defined by
\[
   \bar F_n \;=\; -\log \int_W \prod_{i=1}^n p(X_i\mid w)\,\varphi(w)\,dw .
\]
The integral is the probability the model assigns to the sample before any parameter is chosen: the probability $\prod_i p(X_i\mid w)$ that parameter $w$ assigns, averaged over the prior. So $\bar F_n$ is small when the model, prior included, predicted the data well; and it is random, since it depends on the sample. (This is not the variational free energy $F(q,x)$ of Section~\ref{sec:epistemics}, which scores a belief against a single observation: that one bounds a negative log evidence from above, and this one is a negative log evidence.)

Take first the \emph{regular} case, which imposes two conditions. The model is \emph{identifiable}: distinct parameters give distinct distributions, so the optimal set is a single point $W_0=\{w_0\}$. And the Hessian of $K$ at $w_0$, which is the \emph{Fisher information} $I(w_0)$, is positive-definite. Under both conditions the integrand defining $\bar F_n$ concentrates in a shrinking neighborhood of $w_0$; expanding $\log p(X_i\mid w)$ to second order there makes the integral Gaussian, contributing a volume factor $(2\pi/n)^{d/2}\det I(w_0)^{-1/2}$. Taking $-\log$,
\[
   \bar F_n \;=\; nS_n \;+\; \tfrac{d}{2}\log n \;+\; O_p(1),
   \qquad S_n=-\tfrac1n\textstyle\sum_i\log q_0(X_i),
\]
where $S_n$ is the empirical entropy and the remainder is bounded in probability.\footnote{$O_p(1)$ means \emph{bounded in probability}: for every $\varepsilon>0$ there is an $M$ with $\Pr(|\,\cdot\,|>M)<\varepsilon$ for all large $n$.} The first term is the score the true distribution itself would earn on the sample, and no model can beat it on average. The second is the price of not knowing $w_0$: the posterior concentrates in a ball of radius of order $n^{-1/2}$ about $w_0$, whose prior volume is of order $n^{-d/2}$, and minus the log of that volume is $\tfrac d2\log n$. Every parameter costs half a $\log n$, so this term is the \emph{complexity penalty}: what a model with more parameters pays, when models are compared by free energy, for its extra freedom to fit. Watanabe's theorem is the singular replacement, in which $d/2$ gives way to an invariant of the singular geometry.

That invariant is defined by volume growth. Write $v(\varepsilon)=\int_{\{K\le\varepsilon\}}\varphi(w)\,dw$ for the prior volume of the parameters fitting the truth to within $\varepsilon$. For real-analytic $K$ and a smooth positive prior, this volume has the asymptotic form
\[
   v(\varepsilon)\;\sim\;c\,\varepsilon^{\lambda}\,(\log(1/\varepsilon))^{m-1}\qquad(\varepsilon\to 0)
\]
for a constant $c>0$, a rational number $\lambda>0$, and an integer $m\in\{1,\dots,d\}$. The exponent $\lambda$ is the \emph{learning coefficient} of the model and $m$ is its \emph{multiplicity}. In the regular case the set $\{K\le\varepsilon\}$ is an ellipsoid of radius of order $\sqrt\varepsilon$, so $\lambda=d/2$ and $m=1$; in the toy loss $a^2b^2$ above, $\lambda=1/2$ and $m=2$. Smaller $\lambda$ means near-optimal parameters are more plentiful.

\begin{theorem*}[Watanabe's free-energy asymptotics~\cite{watanabe2009}]
For a singular model satisfying Watanabe's regularity conditions (real-analyticity and a mild variance bound\footnote{The \emph{relatively finite variance} condition: writing $f(x,w)=\log(q_0(x)/p(x\mid w))$, $K(w)=\mathbb{E}_{q_0} f(X,w)$, and $V(w)=\mathbb{E}_{q_0} f(X,w)^2$, one assumes $V(w)\le C\,K(w)$ near the optimal set: the variance of the log-likelihood ratio is controlled by its mean.}), as $n\to\infty$,
\begin{equation}\label{eq:free-energy}
   \bar F_n \;=\; n\,S_n \;+\; \lambda\log n \;-\; (m-1)\log\log n \;+\; O_p(1),
\end{equation}
where $\lambda$ and $m$ are the learning coefficient and multiplicity of the model. They take the places of the $d/2$ and the $1$ of the regular case. The $O_p(1)$ remainder converges in distribution.\footnote{In the regular case the limit is $C-\tfrac12\chi^2_d$ for a constant $C$, the Bayesian shadow of Wilks' theorem, since the log-likelihood-ratio statistic $2\bigl(\sup_w\textstyle\sum_i\log p(X_i\mid w)-\sum_i\log p(X_i\mid w_0)\bigr)\Rightarrow\chi^2_d$. In a singular model the limit is in general non-Gaussian, a functional of a Gaussian process on the manifold obtained by resolving the singularities of $K$.} The same coefficient controls generalization. Let $\hat p_n$ be the average of $p(\cdot\mid w)$ over the \emph{posterior}, the probability density on $W$ proportional to $\prod_i p(X_i\mid w)\,\varphi(w)$ (the prior reweighted by how well each parameter fits the data). Then the \emph{generalization error} $G_n=\mathrm{KL}(q_0\,\|\,\hat p_n)$ obeys $\mathbb{E}[G_n]=\lambda/n+o(1/n)$.
\end{theorem*}

In the regular case the theorem recovers the expansion above; in a singular model $\lambda\le d/2$, often strictly. The expansion localizes: restrict the integral defining $\bar F_n$ to a neighborhood of one optimal parameter, and the expansion \eqref{eq:free-energy} holds with that neighborhood's own coefficient. Among the optimal parameters (all fitting the truth exactly, all weighted by the prior), the posterior share of the most degenerate neighborhoods, those of smallest local $\lambda$, grows with $n$. In this Bayesian setting the posterior concentrates on the solutions that the formula $\mathbb{E}[G_n]=\lambda/n$ marks as best-generalizing. This is the theory's rendering of the slogan that learning prefers ``simpler'' solutions: simplicity means a smaller learning coefficient $\lambda$, not fewer parameters.

The degeneracy behind a small $\lambda$ is largely \emph{structural}: it comes from directions in parameter space along which the model $p(\cdot\mid w)$, and hence the loss, does not change, and these are what pull $\lambda$ below $d/2$. For instance, if a column of one weight matrix is zero, the corresponding row of the matrix before it can be changed freely.

The theory is exact for \emph{Bayesian} learning, not directly for stochastic gradient descent (gradient descent driven by noisy, small-batch estimates of the gradient). Bridging the two is the problem \emph{Opposing staircases} below. The two are closer than they appear: Khan and Rue~\cite{khan2023} derive stochastic gradient descent, along with many other optimizers, as approximations of a single Bayesian update rule.

There are two ways to compute $\lambda$. The first is by hand: algebraic geometers know $\lambda$ as the \emph{real log canonical threshold} of $K$ relative to $\varphi$, and one computes it by resolving the singularities of $K$ (Hironaka~\cite{hironaka1964}), which is feasible for small models. The second is by sampling. Introduce a knob $\beta>0$ and the \emph{tempered posterior}
\[
   p_\beta(w)\;\propto\;\varphi(w)\,e^{-\beta n L_n(w)},
   \qquad L_n(w)=-\tfrac1n\textstyle\sum_i\log p(X_i\mid w),
\]
which interpolates between the prior ($\beta=0$) and concentration on the best-fitting parameters ($\beta\to\infty$). At $\beta=1/\log n$, the average of $nL_n$ over the tempered posterior exceeds its minimum by $\lambda\log n$, to leading order~\cite{watanabe2013}. So $\lambda$ can be \emph{estimated} by sampling the tempered posterior, needing neither the true distribution nor an explicit resolution of singularities. (In statistics, taking $\beta<1$ hedges against a misspecified model: Gr\"unwald's ``Safe Bayesian'' learns $\beta$ from the data~\cite{grunwald2012}, and Alquier and Ridgway prove that tempered posteriors still concentrate~\cite{alquier2020}.)

Sampling a density in millions of dimensions is not a small matter, but the diffusion $dw_t=-\nabla U(w_t)\,dt+\sqrt2\,dB_t$ with $U=\beta nL_n-\log\varphi$ has stationary density proportional to $e^{-U}$, which is $p_\beta$; run in discrete time with small-batch gradients, it is gradient descent with Gaussian noise added at each step, and this is the sampler used in practice~\cite{lau2023}. That practicality is what lets the theory meet networks far too large to analyze by hand.

What one estimates is not the global $\lambda$ but the \emph{local learning coefficient}: the same invariant computed in a neighborhood of the particular solution $w^\star$ that training actually reached, measuring how degenerate \emph{that} solution is~\cite{lau2023}. As a network trains, its estimated local learning coefficient is observed to jump at discrete moments, called \emph{phase transitions}: the solution moves to a qualitatively new, differently-degenerate part of the landscape as a new piece of internal structure forms. Tracking these jumps is the program of \emph{developmental interpretability}: so far on relatively small networks and with preliminary validation, though related estimators from the same theory have already surfaced tens of thousands of candidate structures inside a billion-parameter language model~\cite{murfet2026}.

Two networks can fit the same training data by different internal mechanisms, and two mechanisms that agree on the training distribution can disagree on inputs beyond it. A held-out test drawn from that distribution cannot tell them apart. In the Bayesian setting the geometry decides which mechanism is selected: the most degenerate solution, the one of smallest $\lambda$, receives the most posterior weight. Whether gradient descent selects the same way is the bridge problem \emph{Opposing staircases} below. Nothing in the theory says that the selected solution implements the behavior intended rather than a shortcut that mimics it on the training distribution and fails beyond. Whether small $\lambda$ favors mechanisms that are also interpretable and robust is the question this program must answer for safety; the case that safety turns on such questions is made in~\cite{lehalleur2025}.

\begin{openproblem}[\href{https://github.com/lionellevine/MAIS/blob/main/open-problems/MAIS-O6.md}{\texttt{MAIS-O6}} (Does geometric simplicity force a legible mechanism?)]
Fix an odd integer $p$ and $H\geq 2p-1$. On input $(a,b)\in(\mathbb{Z}/p\mathbb{Z})^2$, consider the one-hidden-layer network
\[
   f_w(a,b)=\sum_{j=1}^H V_j\bigl(u'_j(a)+u''_j(b)\bigr)^2\in\mathbb{R}^p,
\]
trained, with squared-error loss averaged over all $p^2$ inputs, to output the coordinate vector indexed by $a+b$. The $j$th summand is one hidden neuron: it adds two entries from the lookup tables $u'_j,u''_j\colon\mathbb{Z}/p\mathbb{Z}\to\mathbb{R}$, squares the result, and writes the vector $V_j$ to the output. Fourier inversion gives exact fits in which each neuron that contributes nontrivially carries a single frequency: for some $k$, both lookup tables are a constant plus a linear combination of $\cos(2\pi k\,\cdot/p)$ and $\sin(2\pi k\,\cdot/p)$. Up to a change of basis, this is the algorithm found by reverse-engineering trained networks~\cite{nanda2023}.

Is this Fourier structure forced by singular geometry? In a fixed closed ball containing one of these Fourier fits, prove or refute that every exact fit of smallest local learning coefficient has this single-frequency form, apart from neurons whose contribution vanishes identically. The first non-vacuous case is $p=5$ with nine hidden neurons: use its Fourier symmetries to classify the most singular exact fits, or find a non-Fourier one by computer algebra. A proof would connect geometric simplicity to a mechanism one can read from the weights; a counterexample would show that the two notions of simplicity can part ways.
\end{openproblem}

\subsection{Learning in stages}

Learning does not always proceed smoothly; long flat stretches give way to sudden changes.

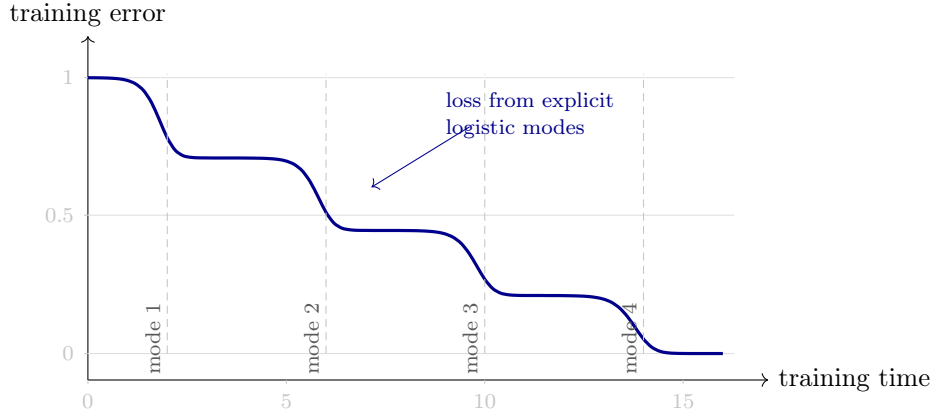
\begin{figure}[ht]
\centering
\begin{tikzpicture}[x=1cm,y=1cm]
  \draw[->] (0,0) -- (9,0) node[right] {\small training time};
  \draw[->] (0,0) -- (0,4.55) node[above] {\small training error};
  \foreach \x/\lab in {0/0,2.625/5,5.25/10,7.875/15} {
    \draw[gray!45] (\x,0) -- (\x,-0.05) node[below,font=\scriptsize] {\lab};
  }
  \foreach \y/\lab in {0.35/0,2.18/0.5,4.0/1} {
    \draw[gray!45] (0,\y) -- (-0.05,\y) node[left,font=\scriptsize] {\lab};
  }
  \draw[gray!25] (0,0.35) -- (8.55,0.35);
  \draw[gray!25] (0,2.18) -- (8.55,2.18);
  \draw[gray!25] (0,4.0) -- (8.55,4.0);
  \draw[blue!55!black,very thick,smooth] plot coordinates {
    (0.000,3.999) (0.105,3.998) (0.210,3.996) (0.315,3.992) (0.420,3.983)
    (0.525,3.962) (0.630,3.918) (0.735,3.831) (0.840,3.673) (0.945,3.443)
    (1.050,3.203) (1.155,3.040) (1.260,2.967) (1.365,2.945) (1.470,2.939)
    (1.575,2.938) (1.680,2.937) (1.785,2.937) (1.890,2.937) (1.995,2.937)
    (2.100,2.936) (2.205,2.935) (2.310,2.933) (2.415,2.928) (2.520,2.918)
    (2.625,2.896) (2.730,2.852) (2.835,2.768) (2.940,2.624) (3.045,2.424)
    (3.150,2.218) (3.255,2.076) (3.360,2.009) (3.465,1.987) (3.570,1.980)
    (3.675,1.979) (3.780,1.978) (3.885,1.978) (3.990,1.978) (4.095,1.978)
    (4.200,1.977) (4.305,1.976) (4.410,1.973) (4.515,1.967) (4.620,1.956)
    (4.725,1.933) (4.830,1.889) (4.935,1.810) (5.040,1.680) (5.145,1.507)
    (5.250,1.333) (5.355,1.210) (5.460,1.149) (5.565,1.127) (5.670,1.120)
    (5.775,1.118) (5.880,1.118) (5.985,1.118) (6.090,1.117) (6.195,1.117)
    (6.300,1.116) (6.405,1.114) (6.510,1.111) (6.615,1.105) (6.720,1.092)
    (6.825,1.069) (6.930,1.026) (7.035,0.951) (7.140,0.836) (7.245,0.688)
    (7.350,0.542) (7.455,0.437) (7.560,0.382) (7.665,0.360) (7.770,0.353)
    (7.875,0.351) (7.980,0.350) (8.085,0.350) (8.190,0.350) (8.295,0.350)
    (8.400,0.350)
  };
  \foreach \x/\name in {1.05/1,3.15/2,5.25/3,7.35/4} {
    \draw[gray!45,densely dashed] (\x,0.35) -- (\x,4.05);
    \node[font=\scriptsize,gray!70!black,rotate=90,anchor=south] at (\x+0.04,0.55) {mode \name};
  }
  \node[font=\scriptsize,align=left,blue!55!black] at (5.85,3.5)
    {loss from explicit\\logistic modes};
  \draw[->,blue!55!black] (5.05,3.35) -- (3.75,2.55);
\end{tikzpicture}
\caption{A learning staircase computed from the deep-linear solution of Saxe et al. In mode $i$, $s_i$ is the target singular value and $a_i(t)$ is the coefficient learned by time $t$, satisfying $\tau\dot a_i=2a_i(s_i-a_i)$. The curve plots normalized training error $\sum_i(s_i-a_i(t))^2$ for $\tau=0.5$ and $s=(1,0.95,0.9,0.85)$, with small initial coefficients chosen so the four modes turn on at separated times. Each mode turning on produces one drop.}
\label{fig:staircase}
\end{figure}

A solvable model of such stepwise learning is the deep \emph{linear} network: a product of matrix factors $W_L\cdots W_1$ fit to input--output pairs $(x,y)$ by \emph{gradient flow} (gradient descent in continuous time) on the squared error $\tfrac12\mathbb{E}\|y-W_L\cdots W_1x\|^2$. Assume the inputs are \emph{whitened} (linearly transformed so that $\mathbb{E}[xx^\top]=I$) and write $\Sigma=\mathbb{E}[yx^\top]$ for the input--output correlation matrix. The system is then
\[
   \tau\,\dot W_\ell \;=\; (W_L\cdots W_{\ell+1})^{\top}\,\bigl(\Sigma-W_L\cdots W_1\bigr)\,(W_{\ell-1}\cdots W_1)^{\top},\qquad \ell=1,\dots,L,
\]
with $\tau$ a time constant and empty products equal to the identity. The composite map is linear, but as a function of the separate factors the loss is a non-convex polynomial, so the flow is nonlinear. Saxe, McClelland, and Ganguli~\cite{saxe2014} show that, for initial weights aligned with the singular vectors of $\Sigma$, the singular value decomposition \emph{decouples} the dynamics into independent scalar equations, one per singular value $s$: for two factors started in balance (equal coefficients in the two factors), each is the logistic equation $\tau\dot a=2a(s-a)$, whose solution rises sigmoidally from near $0$ to $s$ over a window of width of order $\tau/s$, after a delay that grows only logarithmically as the initial value shrinks (Figure~\ref{fig:staircase}). The decoupling and the scalar solutions are exact; the initialization then sets the order of events: starting from small initial weights, modes with larger $s$ turn on earlier, and when their turn-on times are well separated the training error falls not smoothly but in a staircase: long plateaus near saddle points of the loss, each broken by a sharp drop as the next mode is acquired. Here ``which structure the network learns, and in what order'' has a closed-form answer.

The modular-addition task above shows a cousin of the staircase, not the same phenomenon: a network fits its training data early but generalizes only much later, and abruptly. This is \emph{grokking}~\cite{nanda2023}, a delayed jump in performance on \emph{held-out} data after the training loss is already low, whereas the staircase is a stepwise fall of the \emph{training} loss itself. Both are read, in this program, as transitions between differently-singular regions of parameter space.

The staircase survives beyond the linear case: Abbe, Boix-Adser\`a, and Misiakiewicz~\cite{abbe2023} prove that training a two-layer network on a \emph{sparse} target (one depending on only a few of its inputs) again proceeds in stages, saddle to saddle, with the waiting time before each stage set by a combinatorial \emph{leap complexity} measuring how many inputs must be assembled at once. Carrying such exact accounts to the nonlinear networks used in practice is open.

\begin{openproblem}[\href{https://github.com/lionellevine/MAIS/blob/main/open-problems/MAIS-O7.md}{\texttt{MAIS-O7}} (Opposing staircases).]
Watanabe's theorem is Bayesian; gradient descent is not. (What bridging the two would mean for safety is discussed in~\cite{lehalleur2025}.) Build the bridge first in the deep-linear staircase above, where both sides are explicit. Each plateau sits at a saddle where only the $k$ largest modes have been learned, and there the definition itself is part of the problem: the local learning coefficient is defined at a local \emph{minimum}, as the volume exponent of the set of nearby parameters fitting almost as well, but a saddle has descent directions, so that volume does not shrink to zero. Formulate the right local invariant at the $k$th saddle (one candidate, the \emph{two-sided} exponent: the volume of nearby parameters whose loss lies within $\varepsilon$ of the saddle's, above or below, shrinks like $\varepsilon^{\lambda}$), compute it for the deep-linear network, and prove or refute the monotone picture that gives this problem its name: as gradient flow descends the staircase of losses, the invariant that, on the Bayesian side, controls generalization climbs an opposing staircase. As a first check, the estimator of Lau et al.~\cite{lau2023} runs at any parameter, saddle or not; run it along a simulated staircase and see.
\end{openproblem}

\subsection{Generalization beyond the training distribution}

The theory above concerns generalization \emph{in distribution}: test data drawn from the same law as the training data. The safety-critical regime is the opposite one, \emph{out of distribution}: how a trained system behaves on inputs unlike anything it saw in training. Two parameter settings with identical training loss can implement different mechanisms (as above) and therefore extrapolate differently; which extrapolation gradient descent selects is, at present, something we observe after the fact rather than predict.

In \emph{goal misgeneralization}~\cite{langosco2022}, a system under \emph{distribution shift} (a deployment input distribution unlike the one it trained on) keeps its \emph{capabilities} while pursuing a \emph{different goal} from the one intended. The standard example is an agent trained to collect a coin in a video game where the coin always sat at the right end of the level: it learned ``move right,'' not ``reach the coin,'' and the two came apart the moment the coin moved. 

\begin{openproblem}[\href{https://github.com/lionellevine/MAIS/blob/main/open-problems/MAIS-O8.md}{\texttt{MAIS-O8}} (A predictive theory of out-of-distribution generalization).]
In the coin-collecting environment of Langosco et al.~\cite{langosco2022}, reduced to one dimension, two policies fit the training data perfectly: ``move right'' (the proxy) and ``go to the coin'' (the intended goal), which agree whenever the coin sits at the right end. Train a two-layer network by gradient descent on the logistic loss to imitate optimal play, randomizing the coin's position in an $\varepsilon$-fraction of training episodes. Determine the probability, over the random initialization, that the trained network follows the proxy on a probe state where the two policies disagree, as a function of $\varepsilon$, the width, and the input encoding. For linear policies the answer is a theorem---the encoding, not the initialization, decides---and the standard infinite-width limits~\cite{jacot2018,mei2018,chizat2020} follow it under their convergence hypotheses; for networks of finite width it is open at every $\varepsilon$, including $\varepsilon=0$, already at width two.
\end{openproblem}


\section{Getting started in AI safety}
\label{sec:conclusion}

Technological progress is a human choice, not a law of nature. Mathematics holds truth-seeking as a core value, but the civic mathematician also asks what our truths are for, and whether they can help design AI that is more legible, steerable, and cooperative with us.

The open problems in this invitation are meant as entry points into the new field of Math for AI Safety. Many more problems can be found in the \href{https://github.com/lionellevine/MAIS/blob/main/open-problems/README.md}{MAIS repository}, and the \href{https://mbrcic.github.io/ai-safety-formalization-atlas/}{AI Safety Formalization Atlas}, maintained by Mario Br\v{c}i\'c, collects machine-checked Lean proofs of results relevant to AI safety, including verdicts on submitted solutions to MAIS problems. I hope you'll pick a problem whose language already feels like home, strip it down to the simplest interesting case, and start proving things!

The rest of this section collects useful links to the growing ecosystem of AI safety organizations and funders.

\subsection*{Upskilling}
The \href{https://www.arena.education/}{ARENA} curriculum teaches the engineering side; the month-long \href{https://www.iliad.ac/intensive}{Iliad Intensive} is a full-time course on the theory, aimed at mathematicians, physicists, and theoretical computer scientists. \href{https://www.matsprogram.org/}{MATS} and \href{https://sparai.org/}{SPAR} pair newcomers with mentors for a first research project.

\subsection*{Institutes}
Research organizations in AI safety include \href{https://www.alignment.org/}{ARC}, \href{https://humancompatible.ai/}{CHAI}, \href{https://constellation.org/}{Constellation}, \href{https://www.far.ai/}{FAR AI}, \href{https://maisi.org/}{MAISI}, \href{https://resi.org/}{RESI}, and \href{https://resolution.org}{Resolution}. Many run long-term fellowship or visitor programs.

\subsection*{Meetings}
\href{https://www.iliadconference.com}{ILIAD} is a conference on mathematical approaches to AI alignment, and FAR AI runs a series of \href{https://www.far.ai/events}{Alignment Workshops}.

\subsection*{Funding}
Grantmakers supporting work in AI safety include \href{https://coefficientgiving.org/funds/navigating-transformative-ai/}{Coefficient Giving}, \href{https://www.cooperativeai.com/}{CAIF}, \href{https://survivalandflourishing.fund/}{SFF}, and \href{https://alignmentproject.aisi.gov.uk/}{UK AISI}.

\section*{Acknowledgments}

I thank Jesse Hoogland, Hyojeong Son, Jacob Tsimerman, Claude, and Codex for many inspiring conversations. Scott Aaronson, Ahmed Bou-Rabee, Vince Conitzer, Jonathan Gabor, Chris Hillar, Brad Knox, Phil Sosoe, Samuel Speas, Kate Stange, Steve Strogatz, Ariel Yadin, and an anonymous referee provided valuable feedback on an early draft.

\section*{AI collaborators}
\label{sec:ai-collaborators}

I wrote the first draft in collaboration with Claude Opus 4. To start, I instructed Opus to read the full text of approximately 100 AI safety papers and produce a summary of the main results and techniques of each, along with possible relevance to this invitation paper. I supplied the overarching structure of the invitation (organized by field, ending each section with an open problem) and wrote parts of the introduction. Opus then used its summaries along with samples of my writing to write a draft of each section in my voice, which I edited extensively. Later revisions were polished with the help of Claude Fable 5, GPT 5.6 Sol, and GPT 6 Astra.

In July 2026 an exuberant Fable produced approximately 200 pages worth of open problems in a single night while I was asleep. I asked Sol to audit these for openness, well-posedness, and plausible AI safety relevance. The problems that passed this audit were used to seed the \href{https://github.com/lionellevine/MAIS/blob/main/open-problems/README.md}{master list of open problems} in the \href{https://github.com/lionellevine/MAIS}{Math for AI Safety (MAIS) repository}, a new hub for open-source collaboration. 

The audit also turned up a small surprise: while checking the problem that became \texttt{MAIS-O1}, Sol found a gap in the published proof of Critch's bounded L\"ob theorem, and proposed a repair (\href{https://github.com/lionellevine/MAIS/tree/main/papers/P2}{MAIS-P2}). Critch corrected the hypothesis (\href{https://x.com/AndrewCritchPhD/status/2082141000831795467}{post on X, July 28, 2026}) and a version of his theorem, in a bespoke proof calculus, has been proved in Lean~\cite{duclaux2026}.

\appendix

\section{Glossary of machine learning terms for mathematicians}
\label{app:glossary}

\begin{description}
\item[\glsentry{neural-network}{neural network}] a function with many (up to trillions of) adjustable real parameters (its \emph{weights}), built by composing simple layers and tuned by \emph{gradient descent} to reduce a training \emph{loss}.
\item[\glsentry{layer}{layer}] a map $x\mapsto\sigma(Wx+b)$: an affine map, whose matrix and vector entries are weights, followed by a fixed nonlinearity $\sigma$ applied to each coordinate; a network is a composition of layers (a transformer alternates them with attention layers, which mix across positions).
\item[\glsentry{neuron}{neuron}] one coordinate of a layer, the scalar function $x\mapsto\sigma(\langle w_j,x\rangle+b_j)$ read off row $j$ of the layer's matrix; in a transformer, a coordinate of one of its non-attention layers.
\item[\glsentry{language-model}{language model}] a neural network trained to predict the next token of text; today's chat systems are language models further trained to follow instructions.
\item[\glsentry{token}{token}] the atomic unit (word or word-fragment) a language model reads and predicts.
\item[\glsentry{loss}{loss}] the real-valued function a network is trained to minimize; for language models, the negative log-probability it assigned to the actual next token.
\item[\glsentry{weights}{weights}] the trainable parameters of an artificial neural network: the entries of its matrices, trained by gradient descent and fixed at deployment time.
\item[\glsentry{activations}{activations}] the values that flow through a neural network on a given input (the outputs of its layers), as opposed to the fixed weights.
\item[\glsentry{activation-space}{activation space}] the copy of $\mathbb{R}^n$ in which a given layer's activations live, one coordinate per neuron; feature directions are directions in this space.
\item[\glsentry{gradient-descent}{gradient descent}] the optimizer: repeatedly nudge the weights $w$ against the gradient of the loss, $w \leftarrow w - \eta\,\nabla L(w)$.
\item[\glsentry{relu}{ReLU}] the ``rectified linear unit'' $x \mapsto \max(0,x)$, a canonical piecewise-linear nonlinearity; modern transformers often use smooth or gated variants instead (GELU, SiLU, SwiGLU).
\item[\glsentry{transformer}{transformer}] the dominant neural network architecture for language: a stack of layers that mix information across token positions (each position selectively reads from the others) and read from / write to the residual stream.
\item[\glsentry{residual-stream}{residual stream}] the running vector of activations a transformer reads from and writes to at each layer; a natural home for feature directions.
\item[\glsentry{feature}{feature}] a hypothesized variable or property a network represents in its activations (e.g.\ ``the text is in French'').
\item[\glsentry{probe}{probe}] a simple (usually linear) function of a network's activations, trained to read off some quantity of interest.
\item[\glsentry{policy}{policy}] an agent's decision rule (a possibly randomized map from observations to actions).
\item[\glsentry{reinforcement-learning}{reinforcement learning}] training an agent by scoring its actions with a numerical \emph{reward} and adjusting its policy to earn more of it.
\item[\glsentry{superposition}{superposition}] a hypothesized way for a neural network to represent more features than it has dimensions (as non-orthogonal directions, tolerable when few features are active at once).
\item[\glsentry{sparse-autoencoder}{sparse autoencoder}] a network trained to reconstruct another network's activations as sparse combinations of learned dictionary directions; an empirical tool for extracting features.
\item[\glsentry{interpretability}{interpretability}] reverse-engineering a trained network into human-understandable structure (circuits, features, algorithms), validated by intervention and not by human-readable description alone.
\item[\glsentry{chain-of-thought}{chain of thought}] the text a language model writes while reasoning toward its answer; readable by humans, and for now the main window onto an AI's reasoning.
\item[\glsentry{activation-steering}{activation steering}] changing a network's behavior at runtime by adding a fixed vector to its activations, without retraining; see Section~\ref{sec:interp}.
\item[\glsentry{gradual-disempowerment}{gradual disempowerment}] the risk that humans lose influence over the institutions that shape our lives, by incremental delegation to AI systems we do not understand rather than by any sudden loss of control.
\end{description}

\end{document}